\documentclass[10pt,twoside]{article}

\usepackage{fancyhdr}    
\newlength{\myparskip}              \newlength{\myparindent}             %
\newlength{\oldparskip}             \newlength{\oldparindent}            %
\parskip\myparskip                  \parindent\myparindent               %
\fancypagestyle{plain}{%
   \fancyhead[LO]{} 
   \fancyhead[CO]{}\fancyhead[RO]{}
}

\newenvironment{proof}{\topsep=\smallskipamount \partopsep=0pt  %
 \begin{trivlist} \itemindent=\parindent                        %
 \item[\hskip \labelsep\emph{Proof:}]}{\qed\end{trivlist}}      %
\let\qed=\relax                                                 %
\def\qed                                                        %
 {{\unskip\nobreak\hfil\penalty50                               %
   \quad\hbox{}\nobreak\hfil $\Box$                             %
   \parfillskip=0pt \finalhyphendemerits=0 \par}}               %

\newtheorem{lemma}{Lemma}

\newtheorem{theorem}{Theorem}
\newtheorem{remark}{Remark}

\newcommand{\institute}[1]{\newcommand{\where}{#1}}

\renewcommand{\and}{\text{and }}

\title{Calorons, Nahm's equations on $S^1$ and bundles over $\PPP$}
\author{Benoit Charbonneau
   \and Jacques Hurtubise
\thanks{The first author is supported by an NSERC PDF and wishes to
thank Alexandru Ghitza for useful discussions.  The second author is
supported by NSERC and FQRNT grants.  The diagrams in this paper were
created using Paul Taylor's Commutative Diagram package.}}
\institute{Department of Mathematics and Statistics, McGill University\\
           805 Sherbrooke St. W, Montreal, Quebec, H3A 2K6, Canada.\\
           E-mails: benoit@alum.mit.edu and jacques.hurtubise@mcgill.ca}
\date{October 26, 2006}
\usepackage{amssymb,amsmath}
\input diagrams  
     \diagramstyle[size=1.8em]
     \newarrow{Iso}C---{>>}
\usepackage[hyperindex=true,pdfauthor={Benoit Charbonneau (benoit@alum.mit.edu) and Jacques Hurtubise (jacques.hurtubise@mcgill.ca)},pdftitle={Nahm's equations on S1 and bundles over P1xP1},pdfsubject={July 2006}]{hyperref}
\DeclareMathAlphabet{\mathdj}{U}{msb}{m}{n}
\newcommand{\R}{\ensuremath{\mathdj {R}}}  
\newcommand{\C}{\ensuremath{\mathdj {C}}}  
\newcommand{\Z}{\ensuremath{\mathdj {Z}}}  
\newcommand{\N}{\ensuremath{\mathdj {N}}}  
\newcommand{\PP}{\ensuremath{\mathdj {P}}} 
\newcommand{\PPP}{{\PP^1\times\PP^1}}      
\newcommand{\A}{\mathcal{A}}               
\newcommand{\B}{\mathcal{B}}               
\newcommand{\MM}{\mathbf{M}}               
\DeclareMathAlphabet{\Gothique}{U}{euf}{m}{n}
\newcommand{\FF}{{\Gothique F}}            
\newcommand{\OO}{{\mathcal{O}}}            
\newcommand{\vect}[1]{{\begin{bmatrix}#1\end{bmatrix}}} 
\newcommand{\col}[1]{{\begin{matrix}#1\end{matrix}}}    
\newcommand{\Ro}{R^1\pi_*}                 
\renewcommand{\Im}{\mathrm{Im}}            
\newcommand{\coker}{\mathrm{coker}}        
\newcommand{\codim}{\mathrm{codim}}        
\newcommand{\rk}{\mathrm{rk}}              
\newcommand{\Aut}{\mathrm{Aut}}            
\newcommand{\SLt}{Sl(2,\C)}                
\newcommand{\flo}{E_{0,+}}                 
\newcommand{\flinf}{E_{\infty,-}}          
\newcommand{\Eof}{E_{0+}}                  
\newcommand{\Eif}{E_{\infty -}}            
\newcommand{\ra}{{k_1}}                    
\newcommand{\rb}{{k_2}}                    
\newcommand{\rc}{{k_3}}                    
\newcommand{\rd}{{k_4}}                    
\newcommand{\map}[2]{\Phi_{#1,#2}}  
\newcommand{\theproofcomplete}{The proof is now complete.} 
\newcommand{\VE}[1]{V_{#1}}                        
\newcommand{\VK}[1]{W_{#1}}                        
\newcommand{\VKi}[1]{\overline{W}_{#1}}            
\newcommand{\VKK}[1]{\overline{\overline{W}}_{#1}} 

\begin{document}
\maketitle
\begin{abstract}
The moduli space of solutions to Nahm's equations of rank $(k,k+j)$
on the circle, and hence, of $SU(2)$ calorons of charge $(k,j)$, is
shown to be equivalent to the moduli of holomorphic rank 2 bundles
on $\PPP$ trivialized at infinity
 ($\{\infty\}\times\PP^1\cup\PP^1\times\{\infty\}$) with $c_2=k$ and
equipped with a flag of degree $j$ along $\PP^1\times\{0\}$.
 An explicit matrix description of these spaces
is given by a monad construction.
\end{abstract}

\section{Introduction.}
The four-dimensional (anti-)self-dual Yang--Mills (ASD) equations, and their
solutions (called \emph{instantons}),
 are by now a
staple of both geometry and physics, whose myriad uses and properties are too
many to summarise here. The earliest base manifold on which these
equations were studied was simply $\R^4$;
various constructions in this case have been given, the
most efficient of which is the well known ADHM construction \cite{ADHM}.

Early on, solutions were produced by reducing the equations under
the various symmetry groups acting on $\R^4$ (e.g. \cite{HSV,HV1978}).
These reduced equations turn out to be interesting in their own right. Indeed,
invariance under the action of $\R$ by translation produces
monopoles, solutions to the Bogomolny equations (\cite{ward-charge2});
invariance under the action of $\R^2$ produces the Hitchin
equations whose analysis  tells us a
lot about bundles on Riemann surfaces (\cite{HitchinSDRiemann}). Invariance
under $\R^3$ yields Nahm's equations, some important ordinary differential
equations. The case (or rather two cases as we will see)
that concerns us here, that of minimal (translation)
invariance, under a single discrete translation ($\Z$-invariance),
with suitable boundary conditions, corresponds to the case of
\emph{calorons}. These gauge fields have seen a recurrence of
interest recently, for a variety of reasons; see
\cite{Bruckmann2003,Bruckmann2004,Bruckmann2002,Kraan2000,Kraan1998Tduality,kraanBaalInside,kraan1998,lee1998,leelu1998,leeyi2003,Norbury2000}.

{}From this list of examples, it would seem that the most interesting cases
were produced by
considering various Abelian groups acting on $\R^4$ by translation.
While this is a question of taste, these cases do possess a most interesting
feature,
a correspondence due to Nahm (\cite{Nahm,corrigan-goddard}), only proven in
certain cases, which postulates an isomorphism between the moduli of
instantons on  $\R^4$ invariant under the action of a closed subgroup $G$
of translations of $\R^4$  (which can be thought of as suitable fields
on $\R^4/G$), and  instantons on  $\R^4$ invariant under the action of
a dual group  $\hat G$.  The
boundary conditions for both sets of fields must be defined with care for
the correspondence to hold. So
far quite a good set of cases are known; see \cite{jardimsurvey}.

In the case of calorons, the heuristic suggests a correspondence
between calorons (instantons on $\R^4/\Z = S^1\times \R^3 $, with
appropriate boundary conditions), and solutions to {\it Nahm's
equations} (o.d.e.s given by reducing the ASD equations) on the
circle, again with suitable boundary behavior at selected points on
the circle  (see below). This correspondence has been partially
proved by Nye \cite {nye} and Nye--Singer (\cite{singernye}).
We complete the proof in \cite{benoitjacques2}, showing  the
correspondence  is exact.

Thus the moduli space of calorons  corresponds
to the moduli space of solutions to Nahm's equations on the circle,
and it is this space that we examine in this paper. In deciding what
this space should be, we are guided by a few basic ideas.

The first of these is that calorons with gauge group $K$ can be thought as
monopoles over $\R^3$, with values in the  Kac--Moody algebra $\tilde LK$.
This point of view has been developed by Garland and Murray
(\cite {garland-murray}), and is a very useful way of understanding calorons,
in particular for moduli.  Indeed, our second idea is that monopoles
with compact simple gauge group $K$ and maximal symmetry breaking at infinity
(part of the boundary conditions) on $\R^3$
correspond to rational maps of $\PP^1$ to $K/T$, with $T$ the maximal torus.
This is known for $K$ a classical group,
(\cite {DonaldsonSU2,hurtubiseClassification}) but should hold
for  all compact simple groups. The space $K/T$ can also be written as
$K_\C/B$, where $K_\C$ is the complexification of $K$, and $B$ a Borel
subgroup.

Combining these two ideas, the moduli space of calorons for gauge group $K$,
(or rather of solutions to Nahm's equations) should be that of rational maps
into the homogeneous space given by quotienting the loop group $LK_\C$ of
$K_\C$ by the subgroup $\hat LB$ consisting of Fourier series with only terms
of non-negative degree, and with the degree zero term lying in $B$.

This is where a third idea, due to Atiyah (\cite {Atiyah1984}) comes into play.
One thinks of elements of $LK_\C$ as transition
functions for a $K_\C$-bundle on $\PP^1$; a map of $\PP^1$ into $LK_\C$
then defines a bundle
over $\PP^1\times\PP^1$. Working through the quotienting by $\hat LB$, as in
\cite{Atiyah1984}, tells us the theorem that we are going to prove in this
paper. We restrict to the case of $SU(2)$ calorons.

\begin{theorem}\label{MainTheorem}
The moduli space of solutions of rank $(k,k+j)$ to Nahm's equations on
the circle (Equations (\ref{eqn:Nahmcomplex}) and
(\ref{eqn:Nahmreal})) is equivalent to the moduli space of pairs
vector bundles $E$ of rank 2 on $\PP^1\times\PP^1$, first Chern
class zero and second Chern number $k$, trivialized over
$\{\infty\}\times \PP^1 \cup \PP^1\times \{\infty\}$, and equipped
with an injective map $\phi$ from the line bundle ${\OO}(-j)$ of
degree $-j$ to the restriction of $E$ to $ \PP^1\times\{0\}$
(up to non-zero scalar multiple), such
that the image of $\phi(\infty)$ lies in the subspace spanned by the
second vector of the trivialization at $(\infty,\infty)$.
\end{theorem}

The proof of this theorem goes in several steps, each interesting in its own
right:
\begin{itemize}
\item We first show that the pairs $(E,\phi)$ are equivalent to certain
matrices, satisfying some algebraic conditions. This step is a
generalization of the monad construction of \cite{Donaldson1984}. It
is the subject of Sections \ref{sec:monadconstruction} and \ref{sec:proofMonads}.
\item We then show that the monad is equivalent to a set of sheaves
on $\PP^1$, and maps between them. We do so in Section \ref{sec:seqsheaves}.
\item Finally, in
Section \ref{sec:Nahmcomplex}, we show that these sheaves are
equivalent to a Nahm complex over the circle, and hence to a
solution of the Nahm equations.
\end{itemize}

Of the various steps in the chain of equivalences, perhaps the one expressing
the caloron, or bundle plus map, as a diagram of sheaves over $\PP^1$ is most
deserving of comment. The twistor construction of calorons, at least from a
Kac--Moody point of view, is given by some algebro-geometric data over
$T\PP^1$, and the diagram of sheaves is given by restricting this data
to a fiber $\C $ of the projection $T\PP^1\rightarrow\PP^1$, and extending
it to $\PP^1$. This is quite similar to what happens for monopoles

One sees again the theme of Kac--Moody groups.  Our moduli space, let us not
forget, is supposed to correspond to that of maps from $\PP^1$ to a
homogeneous space corresponding to this group. In the finite dimensional
case of homogeneous spaces for $Gl(n,\C)$, one can describe maps into a
flag manifold in terms of similar diagrams of sheaves
(\cite{hurtubiseClassification}); it is not surprising that this pattern
reoccurs here, and that it is a case in which the ``finite-dimensional''
aspects of the Kac--Moody group (root spaces, etc) predominate.

\section{Monad Construction.}\label{sec:monadconstruction}
Let us use standard affine coordinates $(x,y)$ on $\PPP$,
denote $\pi$ the projection on the first factor, and $i_{y_0}\colon \PP^1
\to \PPP$  the injection $x\mapsto (x,y_0)$.  Set $H_1:=\{\infty\}\times\PP^1$
and $H_2:=\PP^1\times\{\infty\}$.
For any
sheaf $\FF$ on $\PPP$, denote
\[\FF_{y_0}:=i_{y_0}{}_*i_{y_0}^*\FF\]
the extension by zero of the restriction of $\FF$ at level $y_0$ and
\[\FF(p,q):=\FF\otimes \OO(pH_1+qH_2).\]

Let $E$ be a $\SLt$-bundle over $\PPP$, with $c_2(E)=k$, trivial over
the fiber $\{\infty\}\times \PP^1$, trivialized over the section
$\PP^1\times\{\infty\}$ (thus equipped with a standard degree zero
flag $\Eif\subset E_\infty$ defined by the first basis vector)
and with given flag $\Eof\subset E_0$ of degree $j$ over $\PP^1\times\{0\}$,
and such that identifying the fiber of $E$ at $(\infty,0)$ and
$(\infty,\infty)$, the flags $\Eof$ and $\Eif$  are transverse.

We define three locally free sheaves $K_0,K_\infty$ and $K_{0\infty}$ by the
exact sequences
\begin{align}
0\rTo K_{0\phantom{\infty}}\rTo  E \rTo &E_0/\flo \cong \OO(j,0)_0  \rTo 0,\label{seq:defK0}\\
0\rTo K_{\phantom{0}\infty} \rTo   E \rTo &E_\infty/\flinf \cong \OO_\infty \rTo 0,\label{seq:defKi}\\
0\rTo K_{0\infty} \rTo  E \rTo &E_0/\flo\oplus E_\infty/\flinf  \rTo 0.\label{seq:defK0i}
 \end{align}
As a consequence we have supplementary sequences
\begin{align}
0\rTo K_0\ \rTo  K_{0\infty}(0,1) \rTo &\OO_\infty  \rTo 0,\label{seq:defK3a}\\
0\rTo K_\infty\rTo  K_{0\infty}(0,1) \rTo &\OO(-j,0)_0   \rTo 0.\label{seq:defK3b}
 \end{align}

In Section \ref{sec:Nahmcomplex}, we define the Nahm complex using
the diagram of sheaves
\begin{equation}\label{Nahm-cplx-diagram}
\begin{diagram}[size=1.5em]
E(0,-1) &\lTo&K_\infty (0,-1)\\
\uTo&&\dTo\\
K_0(0,-1)&\rTo&K_{0\infty},
\end{diagram}
\end{equation}
taking direct images onto $\PP^1$. We have the easily proven lemma.
\begin{lemma} \label{lemma-on-direct-images}
1) $H^0(\PP^1,R^1\pi_*(E(0,-1)) = \C^k$, and $R^1\pi_*(E(0,-1))$ is supported
on $k$ points, counted with multiplicity.

2) $H^0(\PP^1,R^1\pi_*(K_{0\infty})) = \C^{k+j}$, and $R^1\pi_*(K_{0\infty})$ is
supported on $k+j$ points, counted with multiplicity.

3) $R^1\pi_*(K_{\infty}(0,-1))$ and $R^1\pi_*(K_{0}(0,-1))$ are supported
over the whole line; generically, $R^1\pi_*(K_{\infty}(0,-1)) = \OO(k)$,
$R^1\pi_*(K_{0}(0,-1))= \OO(k+j)$.

4) These sheaves fit in the exact sequences
\begin{equation}\label{direct-image-sequences}
\col{
0\rightarrow& \OO(j)&\rightarrow& R^1\pi_*(K_{0}(0,-1))&\rightarrow& R^1\pi_*(E(0,-1))&\rightarrow&0,\\
0\rightarrow& \OO&\rightarrow& R^1\pi_*(K_{\infty}(0,-1))&\rightarrow& R^1\pi_*(E(0,-1))&\rightarrow&0,\\
0\rightarrow& \OO&\rightarrow& R^1\pi_*(K_{0}(0,-1))&\rightarrow& R^1\pi_*(K_{0\infty})&\rightarrow&0,\\
0\rightarrow& \OO(-j)&\rightarrow& R^1\pi_*(K_{\infty}(0,-1))&\rightarrow& R^1\pi_*(K_{0\infty})&\rightarrow&0.
}
\end{equation}
In particular, the natural maps from $R^1\pi_*(K_{\infty}(0,-1))$
and $R^1\pi_*(K_{0}(0,-1))$ to  $R^1\pi_*(E(0,-1))$ and
$R^1\pi_*(K_{0\infty})$ are all surjections, and the kernels are all torsion
free.
\end{lemma}

\begin{proof}
On a generic fiber of $\pi$, we have
$E = \OO\oplus\OO$, $K_0=\OO\oplus\OO(-1)$, $K_\infty=\OO\oplus\OO(-1)$, and
$K_{0\infty}=\OO(-1)\oplus\OO(-1)$.  Hence for
$F= E(0,-1)$, $K_0(0,-1)$, $K_\infty(0,-1)$, or $K_{0\infty}$,
and for $w$ generic, $F|_{\pi^{-1}(w)}$ has  no global sections.  Thus
$\pi_*F=0$,
and $\Ro E(0,-1)$, $\Ro K_{0\infty}$ are torsion, while the sheaves
$\Ro K_{0}(0,-1)$, and
$\Ro K_{\infty}(0,-1)$ are line bundles over the generic
set of $w$ for which $K_0=\OO\oplus\OO(-1)$,  and
$K_\infty=\OO\oplus\OO(-1)$; for
the generic bundle, this is all of $\PP^1$.

Then statement 4 follows from taking the direct image of the sequences
(\ref{seq:defK0})--(\ref{seq:defK0i}).  Statements 1, 2, and 3 follow from the
Grothendieck--Riemann--Roch theorem, and Sequences (\ref{seq:defK0})--(\ref{seq:defK0i}).
\theproofcomplete
\end{proof}

Let $F$ be a vector bundle on $\PPP$ of rank $r$ with first and second Chern
classes
$c_1$ and $c_2$.  Using the Riemann--Roch theorem,
we find
\begin{gather}
\chi(F)=\frac12c_1^2-c_2+(H_1+H_2).c_1+r,\label{chiE}\\
\chi(F(p,q))=\chi(F)+r(p+q)+pH_1.c_1+qH_2.c_1+rpq.\label{chiEpq}
\end{gather}

\begin{lemma}\label{lemma:constant}
Let $F$ be a vector bundle on $\PPP$, trivial on a section
$\PP^1\times\{y\}$.  Then  $H_2.c_1=0$, and as functions of $q$,
\begin{align*} &h^0(F(p,q)) \text{ is constant  for fixed } p\leq-1,\\
               &h^1(F(-1,q)) \text{ is constant,}\\
           &h^2(F(p,q))  \text{ is constant for fixed } p\geq-1.
\end{align*}
Similarly, for  $F$ trivial on a fiber $\{x\}\times \PP^1$,
we have $H_1.c_1=0$, and
\begin{align*} &h^0(F(p,q)) \text{ is constant for fixed } q\leq-1,\\
              &h^1(F(p,-1))  \text{ is constant as functions of $p$,}\\
           &h^2(F(p,q))  \text{ is constant for fixed } q\geq-1.
\end{align*}
\end{lemma}

\begin{proof}
Suppose $F$ is a rank $r$ vector bundle on $\PPP$ which is trivial on
$ \PP^1\times \{y\}$.
We have an exact sequence of sheaves
$0\to F(0,-1) \to F \to F_y \to 0$.
Tensored by $\OO(p,q)$, this sequence gives us a long exact sequence
in cohomology, which ensures
an isomorphism between $H^i(F(p,q))$ and $H^i(F(p,q-1))$ whenever both
$H^i(F_y(p,q))$ and $H^{i-1}(F_y(p,q))$ are $0$. Because of the equality
$H^i(F_y(p,q))=H^i(\PP^1,\OO(p)^r)$, this condition happens precisely to
ensure the $h^i$ are constant as specified.

Since $F$ is trivial along a section, we just proved that
\[\chi(F(-1,q))=\chi(F)-r-H_1.c_1+qH_2.c_1\]
is constant as a function of $q$.  Hence $H_2.c_1=0$.  With this remark, the
first half of this lemma is proved.  The second half is
obtained by symmetry.
\end{proof}

Our bundle $E$, being trivial over the fiber and section at
$\infty$, has thus $c_1=0$.

\begin{lemma}\label{lemma:chernchi}
For  $K_0$, $K_\infty$ and $K_{0\infty}$ defined by
Sequences (\ref{seq:defK0}), (\ref{seq:defKi}), and (\ref{seq:defK0i}),
\begin{align*}  c(K_0)       &=1-H_2+(k+j)H_1H_2,\\
                c(K_\infty)  &=1-H_2+kH_1H_2,\\
                c(K_{0\infty})&=1+2H_2+(k+j)H_1H_2 \end{align*}
and
\begin{align*}\chi(K_0(p,q))       &=-(k+j)+(1+p)(1+2q),\\
              \chi(K_\infty(p,q))   &=-k+(1+p)(1+2q),\\
              \chi(K_{0\infty}(p,q))&=-(k+j)+2q(1+p).\end{align*}
\end{lemma}

\begin{proof}On the following exact sequences, use that $c(\FF_3)=c(\FF_1)c(\FF_2)^{-1}$ when $0\to \FF_1\to \FF_2\to \FF_3\to0$ is
exact:
\begin{gather*}
0\to E(0,-1)\to E\to E_0\to 0,\\
0\to \OO(-j,-1)\to \OO(-j,0)\to \OO(-j,0)_0\to 0, \\
0\to\OO(-j,0)_0\to E_0\to E_0/E_{0+}\to 0,\\
0\to K_0\to E\to E_0/E_{0+}\to 0,\\
0\to K_{0\infty}\to K_0\to E_\infty/E_{\infty-}\to 0. \end{gather*}

Using Equations (\ref{chiE}) and (\ref{chiEpq}), we obtain from
those Chern classes the Euler characteristics.  Setting $j=0$ gives
the answers for $K_\infty$.
\end{proof}

Using this result, we can compute many of the cohomology groups of the
bundles $E, K_0, K_\infty, K_{0\infty}$.

\begin{theorem}[Vanishing] \label{thm:vanishing}
The cohomology groups of the bundles
$E$, $K_0$, $K_\infty$, and $K_{0\infty}$ defined by Sequences (\ref{seq:defK0}), (\ref{seq:defKi}), and (\ref{seq:defK0i}) vanish
as follows:
\begin{gather}
p\leq-1\text{ or }q\leq -1                         \label{eqn:vanishingEK00}
     \kern-4pt\implies\kern-4pt h^0(E(p,q))=h^0(K_0(p,q))=h^0(K_\infty(p,q))=0,\\
p\leq -1\text{ or }q\leq 0                         \label{eqn:vanishingK0i0}
      \implies h^0(K_{0\infty}(p,q))=0,\\
p\geq -1\text{ or }q\geq -1                         \label{eqn:vanishingE2}
                         \implies h^2(E(p,q))=0,\\
p\geq -1\text{ or }q\geq 0\kern-4pt\implies \kern-4pt
h^2(K_0(p,q))=h^2(K_\infty(p,q))=h^2(K_\infty(p,q))=0.
\end{gather}
When $h^0=h^2=0$, we get an exact formula for $h^1$:
\begin{gather}
                                                   \label{eqn:vanishingE1}
\left.\begin{matrix}(p\leq -1\text{ and }q\geq -1)\\
\text{or }(p\geq-1\text{ and }q\leq -1)\end{matrix}\right\}
                           \implies h^1(E(p,q))=k-2(1+q)(1+p),\\
                                                 \label{eqn:vanishingK01}
\left.\begin{matrix}(p\leq-1\text{ and }q\geq 0)\\
\text{or }\\(p\geq -1\text{ and }q\leq -1)\end{matrix}\right\}
                           \kern-6pt\implies\kern-4pt
                              h^1(K_0(p,q))=(k+j)-(1+2q)(1+p),\\
                                                \label{eqn:vanishingK0i1}
\left.\begin{matrix}(p\leq-1\text{ and }q\geq 0)\\
\text{or }(p\geq -1\text{ and }q\leq 0)\end{matrix}\right\}
                           \implies h^1(K_{0\infty}(p,q))=(k+j)-2q(1+p).
\end{gather}
Equation (\ref{eqn:vanishingK01}) is valid for $K_\infty$ by setting $j=0$.
When $j\geq 1$, we get extra information for $K_0$: for $p\leq j-1$, we have
\begin{gather}
\label{eqn:vanishingK00j} h^0(K_0(p,0))=0\\
\label{eqn:vanishingK01j} h^1(K_0(p,0))=k+j-1-p.
\end{gather}
\end{theorem}

\begin{proof}
Lemma \ref{lemma:constant} tell us that $h^2(E(p,q))$ is constant
in the region \[\{(p,q)\mid p\geq -1\text{ or }q\geq -1 \}.\]  For any $i>0$,
and $N$
big enough, Theorem B of Serre (see \cite[p.700]{gharris}) says
$h^i(F(N,N))=0$.  Thus $h^2(E(p,q))=0$ throughout this region,
proving (\ref{eqn:vanishingE2}).  By Serre duality, we have the corresponding
result for $h^0$.

Restricted to a generic section, $K_0$, $K_\infty$ and $K_{0\infty}$ are
trivial.  Hence the first part of Lemma \ref{lemma:constant} applies.
Restricted to the fiber above $\infty$ however, $K_0$ and $K_\infty$
are $\OO\oplus\OO(-1)$ while $K_{0\infty}$ is $\OO(-1)\oplus\OO(-1)$.

Working as in the proof of Lemma \ref{lemma:constant}, we have for
$F\in\{K_0, K_\infty\}$ an exact sequence
\[\cdots\to H^i(F(p-1,q))\to H^i(F(p,q)) \to H^i(\PP^1,\OO(q)\oplus \OO(q-1))
   \cdots\]
ensuring that
\begin{align*} &h^0(F(p,q)) \text{ is constant  for fixed } q\leq-1,\\
           &h^2(F(p,q))  \text{ is constant for fixed } q\geq0.
\end{align*}
Similarly, we obtain
\begin{align*} &h^0(K_{0\infty}(p,q)) \text{ is constant  for fixed } q\leq0,\\
           &h^2(K_{0\infty}(p,q))  \text{ is constant for fixed } q\geq0.
\end{align*}
Again, we deduct from Theorem B the wanted vanishing for
the $h^0$ and $h^2$.

Whenever $h^0=h^2=0$, we have $h^1=-\chi$.  The exact formula
for $h^1$ thus follow from the Riemann--Roch Equation (\ref{chiEpq})
used with Lemma \ref{lemma:chernchi}.
\end{proof}

Because of all this vanishing, the following theorem
guarantees a monad description of those bundles.

\begin{theorem}[Buchdahl's Beilinson's theorem]    \label{buchdahl}
For any ho\-lo\-mor\-phic vector bun\-dle $F$ on $\PPP$,
the\-re is a spectral sequence
with $E_1$-term
\[\begin{matrix}\OO(-1,-1)^{h^2(F(-1,-1))}&
         \ \ \OO(-1,0)^{h^2(F(-1,0))}\oplus \OO(0,-1)^{h^2(F(0,-1))}\ \
         & \OO^{h^2(F)}\\
\OO(-1,-1)^{h^1(F(-1,-1))}&
         \OO(-1,0)^{h^1(F(-1,0))}\oplus \OO(0,-1)^{h^1(F(0,-1))}
         & \OO^{h^1(F)}\\
\OO(-1,-1)^{h^0(F(-1,-1))}&
         \OO(-1,0)^{h^0(F(-1,0))}\oplus \OO(0,-1)^{h^0(F(0,-1))}
         & \OO^{h^0(F)}
\end{matrix}\]
and  \[E_1^{p,q}\Rightarrow E^{p+q}_\infty
            = \begin{cases}  F,&\text{ if }p+q=0,\\
                             0,&\text{ otherwise.}\end{cases}\]
\end{theorem}

\begin{proof}See \cite[p. 144]{Buchdahl1987}.  Use
\cite[Ex. 6.1, p. 237]{Hartshorne1977} and
\cite[Prop. 6.3, p. 234]{Hartshorne1977} to see that the only
extension $0\to \OO(1,0)(1,0)'\to R\to \OO(0,1)(0,1)'\to0$ on $H_0=\PPP$
is the trivial one.
 Hence we replace the term
$E_1^{-1,q}$ in the $H_n$-analogue of Beilinson's theorem by the direct
sum as in the statement above.\end{proof}

Now let us exploit the vanishing Theorem together with the machinery
of monads.  First, let us check existence.

\begin{theorem}[Existence of monad]\label{prop:monadEgen}
Let $E$ be one of the bundles we are considering, that is trivial on $\{\infty\}\times\PP^1$ and
 on $\PP^1\times\{\infty\}$. The bundles $E, K_0$,$K_\infty$  and $K_{0\infty}(0,1)$ are the cohomology of monads of respective type
\begin{equation}\label{seq:monadEgen}\begin{aligned}
\OO(-1,0)^{k\phantom{+j}}\to \OO(-1,1)^{k\phantom{+j}}&\oplus\OO^{k+2\phantom{+j}} \to \OO(0,1)^k,\\
\OO(-1,0)^{k+j}\to \OO(-1,1)^{k+j}&\oplus\OO^{k+j+3} \to \OO(0,1)^{k+j+1},\\
\OO(-1,0)^{k\phantom{+j}}\to \OO(-1,1)^{k\phantom{+j}}&\oplus\OO^{k+3\phantom{+j}} \to \OO(0,1)^{k+1},\\
\OO(-1,0)^{k+j}\to \OO(-1,1)^{k+j}&\oplus\OO^{k+j+2} \to \OO(0,1)^{k+j}.
\end{aligned}
\end{equation}
\end{theorem}

Before starting the proof, set
\begin{equation}\label{eqn:notationV}
\begin{aligned} V_1(F)&:=H^1(F(-1,-2)),
          &\quad &V_2(F):=H^1(F(-1,-1)),\\
                  V_3(F)&:=H^1(F(0,-2)),
          &      &V_4(F):=H^1(F(0,-1)),\\ \end{aligned}
\end{equation}

\begin{proof}We use extensively the Vanishing Theorem \ref{thm:vanishing}.
Suppose $F\in \{E,K_0,K_\infty\}$.  The cohomology groups appearing in the
$E_1$-term for $F(0,-1)$ are $H^i(F(p,q))$ with $p\in\{-1,0\}$ and
$q\in\{-2,-1\}$.  Then $q\leq -1$ implies those groups with $i=0$ are trivial,
and $p\geq-1$ implies those groups with $i=2$ are also trivial.  For the $K_{0\infty}(p,q)$ appearing, with $p,q\in\{-1,0\}$, the fact that $q\leq 0$ ensures
$h^0=0$ while $p\geq -1$ ensures $h^2=0$.

The
$E_1$-term for $F(0,-1)$ with $F\in\{E,K_0,K_\infty,K_{0\infty}(0,1)\}$ thus reduces to the middle row, so that
$F(0,-1)$ is the cohomology of a monad of the form
\[\OO(-1,-1)\otimes V_1(F)\to \OO(-1,0)\otimes V_2(F)\oplus\OO(0,-1)\otimes V_3(F) \to \OO  \otimes V_4(F).\]
  The exact formula for the  dimensions of those spaces are given in
Theorem \ref{thm:vanishing}, and correspond to the values given in Equation
(\ref{seq:monadEgen}).
Tensoring by $\OO(0,1)$ concludes the proof.\end{proof}


The following lemma can be used to double check the result.
\begin{lemma}\label{lemma:rankchern}For the cohomology $F$ of a monad
\begin{equation}\label{eqn:monadtype}
\begin{diagram} \OO(-1,0)^{\ra}&\rTo &\OO(-1,1)^{\rb}\oplus\OO^{\rc} &\rTo &\OO(0,1)^{\rd},\end{diagram}\end{equation}
we have
\begin{align*}
\rk F &= \rb+\rc-\ra-\rd,\\
c(F) &= 1+(\ra-\rb)H_1+(\rb-\rd)H_2\\
&\phantom{=
1}+\bigl(\rb(\ra+\rd+1-\rb)-\ra\rd\bigr)H_1H_2.\end{align*}
\end{lemma}

\begin{proof}From \cite[Lemma 3.1.2, p. 240]{okonek}, we have
\[c(F)=c\bigl(\OO(-1,1)^\rb\oplus \OO^\rc\bigr)
                     c\bigl(\OO(-1,0)^\ra)^{-1}c\bigl(\OO(0,1)^\rd)^{-1}.\]
Since $c(\OO(p,q)^r)= 1+r(pH_1+qH_2)+r(r-1)pqH_1H_2$, we have
the proof.\end{proof}

\begin{lemma}\label{lemma:restraints}
For the monads of Theorem \ref{prop:monadEgen}, written as
\begin{diagram}\MM(F)\colon \col{V_1\\ \otimes\\ \OO(-1,0)}
&\rTo_{\vect{\alpha_1y+\alpha_0\\\beta_1x+\beta_0}}
&\col{V_2\otimes \OO(-1,0)\\ \oplus \\ V_3\otimes \OO}
&\rTo_{\vect{\mu_1x+\mu_0&\nu_1y+\nu_0}}& \col{V_4\\ \otimes\\ \OO(0,1)},
\end{diagram}
the map $\alpha_1$ is an isomorphism, $\beta_1$ and $\mu_1$ are injective,
 $\nu_1$ is surjective, and $\ker(\nu_1)\cap \Im(\beta_1)=0$.
\end{lemma}

\begin{proof}
Note first that for $F\in\{E,K_0,K_\infty,K_{0\infty}(0,1)\}$, we have
$F|_{\PP^1\times\{\infty\}}$ trivial and $F|_{\{\infty\}\times\PP^1}$
isomorphic to $\OO^2$, $\OO\oplus\OO(-1)$, or $\OO(-1)^2$.

Let us first study the restriction to $\PP^1\times\{\infty\}$.  Consider the
monad
\begin{diagram}\MM(F|)\colon V_1\otimes\OO(-1)
&\rTo_{\vect{\alpha_1\\\beta_1x+\beta_0}}
&\col{V_2\otimes \OO(-1)\\ \oplus \\ V_3\otimes \OO}
&\rTo_{\vect{\mu_1x+\mu_0&\nu_1}}& V_4\otimes \OO.
\end{diagram}
Since $0\to V_1\otimes\OO(-1)\to\ker\to F|\to 0$ is exact,
we have  $H^1(\ker)=H^1(F|)=0$.  
Since
$0\to\ker\to V_2\otimes \OO(-1) \oplus V_3\otimes \OO \to V_4\otimes \OO\to 0$
is exact,
we obtain
\begin{diagram}V_3\otimes H^0(\OO)
&\rTo^{\nu_1}& V_4\otimes H^0(\OO)&\rTo &H^1(\ker)=0,
\end{diagram}
whence $\nu_1$ is surjective.

Consider the monad
\begin{diagram}\MM(F|^*(-1))\colon V_4^*\otimes\OO(-1)
&\rTo
&\col{V_2^*\otimes \OO\\ \oplus \\ V_3^*\otimes \OO(-1)}
&\rTo& V_1^*\otimes \OO.
\end{diagram}
Note that as before, $H^1(\ker)=H^1(F|^*(-1))=0$, but here $H^0=0$ as well.
Thus
\begin{diagram}0=H^0(\ker)&\rTo&V_2^*\otimes H^0(\OO)
&\rTo^{\alpha_1^*}& V_1^*\otimes H^0(\OO)&\rTo &H^1(\ker)=0,
\end{diagram}
being exact, we have that $\alpha_1$ is an isomorphism.

Now let us study the restriction to $\{\infty\}\times\PP^1$.  Consider the
monad
\begin{diagram}\MM(F|(-1))\colon V_1\otimes\OO(-1)
&\rTo_{\vect{\alpha_1y+\alpha_0\\\beta_1}}
&\col{V_2\otimes \OO\\ \oplus \\ V_3\otimes \OO(-1)}
&\rTo_{\vect{\mu_1&\nu_1y+\nu_0}}& V_4\otimes \OO.
\end{diagram}
As before, we have an exact sequence
\begin{diagram}0=H^0(\ker)&\rTo&V_2\otimes H^0(\OO)
&\rTo^{\mu_1}& V_4\otimes H^0(\OO),
\end{diagram}
whence $\mu_1$ is injective.

Consider now
\begin{diagram}\MM(F|^*)\colon V_4^*\otimes\OO(-1)
&\rTo
&\col{V_2^*\otimes \OO(-1)\\ \oplus \\ V_3^*\otimes \OO}
&\rTo& V_1^*\otimes \OO.
\end{diagram}
{}From this exact sequence we get in cohomology the exact sequence
\begin{diagram}V_3^*\otimes H^0(\OO)
&\rTo^{\beta_1^*}& V_1^*\otimes H^0(\OO)&\rTo &H^1(\ker)=0,
\end{diagram}
whence $\beta_1$ is injective.

The monad equation $\mu_1\alpha_1+\nu_1\beta_1=0$ and the injectivity of $\mu_1\alpha_1$ imposes $\ker(\nu_1)\cap \Im(\beta_1)=0$.  The proof of Lemma \ref{lemma:restraints} is now complete.
\end{proof}

Since the construction of Buchdahl is natural, a map $\phi\colon
F\to F'$ induces maps in cohomology $\phi_*\colon V_i(F)\to
V_i(F')$, which we denote $\map i{F,F'}$, or $\Phi_i$ when there is
no risk of confusion.

Using Sequences  (\ref{seq:defK0}), (\ref{seq:defKi}), and
(\ref{seq:defK3a}), twisted by
$\OO(p,q)$,
we have that all the $\map i{K_0,E}$, $\map i{K_\infty,E}$  and
$\map i{K_0,K_{0\infty}(0,1)}$ induced by the injections are surjective.
Moreover, because  $k_1=k_2$ sometimes, four of
those maps are isomorphisms.

Sequence (\ref{seq:defK3b}), tensored by $\OO(p,q)$,
yields the exact sequence
\[H^0(\PP^1,\OO(p-j))\to H^1(K_\infty(p,q))\to H^1(K_{0\infty}(p,q+1))
\to H^1(\PP^1,\OO(p-j)).\]

Depending on $p$ and $j$, we obtain surjectivity and/or injectivity.
We can summarize this information with the diagrams of the following lemma.

\begin{lemma}\label{induced}
The four maps of Diagram (\ref{Nahm-cplx-diagram}) induce monad maps with
the following surjectivity/injectivity properties:
\begin{diagram}
 V_1(K_0)&\rIso &V_1(K_{0\infty}(0,1)) &\quad\quad     & V_3(K_0)&\rOnto&V_3(K_{0\infty}(0,1))\\
  \dOnto &      &    \uInto            &               & \dOnto  &      &\uTo_{\begin{matrix}\text{\small surj. if
                                                               $j\leq 1$,}\\ \text{\small inj. if $j\geq 1$ }\end{matrix}}        \\
V_1(E)   &\lIso& V_1(K_\infty)        &               &V_3(E)  &\lOnto& V_3(K_\infty)\\
\\
 V_2(K_0)&\rIso & V_2(K_{0\infty}(0,1))&               &V_4(K_0)&\rOnto& V_4(K_{0\infty}(0,1))\\
  \dOnto &      &    \uInto            &               & \dOnto &      &\uTo_{\begin{matrix}\text{\small surj. if
                                                               $j\leq 1$,}\\ \text{\small inj. if $j\geq 1$ }\end{matrix}}     \\
V_2(E)   &\lIso& V_2(K_\infty)        &               &V_4(E)
&\lOnto&V_4(K_\infty).
\end{diagram}
\end{lemma}

\begin{remark}\label{remark:symmetries}
The monads given by Theorem \ref{prop:monadEgen}, of the type given by
Equation (\ref{eqn:monadtype}), are uniquely determined
up to the action of \[\Aut\bigl(\OO(-1,0)^\ra\bigr)\times \Aut\bigl(\OO(-1,1)^\rb\oplus \OO^\rc\bigr)\times\Aut\bigl(\OO(0,1)^\rd\bigr);\] see
\cite[Lemma 4.1.3 on p. 276]{okonek}.  This group is exactly
\[Gl(\ra,\C)\times Gl(\rb,\C)\times Gl(\rc,\C)\times Gl(\rd,\C).\]
\end{remark}
Let us exploit those symmetries to give normal forms for the monads of
Theorem \ref{prop:monadEgen}.  Before doing so, set
\begin{equation}\label{def:seplus}
s:=\vect{0&0&\cdots &0&0\\ 1&0 & \cdots &0&0\\ 0&1 &\cdots &0&0\\
            \vdots&\vdots & \ddots& \vdots &\vdots\\
            0&0&\cdots&1&0}\quad\text{ and }\quad e_{+}:=\vect{0&\cdots &0 &1}.
\end{equation}

\begin{theorem}\label{thm:monadE}
For $j>0$, the bundle $E$, $K_0$, $K_\infty$ and $K_{0\infty}(0,1)$ are respectively
the cohomology of the monads
\begin{equation}\label{seq:monadE}\begin{diagram} 
  \OO(-1,0)^k&\rTo_{\vect{A-y\\
                          B-x\\
                           D}}^{\A(E)}
   &\col{\OO(-1,1)^k\\ \oplus\\ \OO^{k+2}}
   &\rTo_{\vect{x-B&A-y& C}}^{\B(E)} &\OO(0,1)^k,
\end{diagram}\end{equation}                        
\begin{equation}\label{seq:monadK0}\raise35pt\hbox{\begin{diagram} 
   \col{\OO(-1,0)^{k+j}\\ \phantom{,}}\kern-17pt &\rTo_{\vect{A-y         &  0\\
                                A'          & -y\\
                                B-x         & 0 \\
                                \vect{B'\\0}&s-x\\
                                D           &0  \\
                                0           &e_+}}^{\A(K_0)}
    &\kern-17pt\col{\OO(-1,1)^{k+j}\\ \oplus \\ \OO^{k+j+3}}\kern -17pt
     &\rTo_{\vect{x-B           & 0   &A-y &0  &C           &0\\
                  \vect{-B'\\0} &x-s  & A' &-y &C'          & 0\\
                   0            &-e_+ & 0  &0  &\vect{1\ 0} & -y}}^{\B(K_0)}
    &\kern-4pt O(0,1)^{k+j+1},
\end{diagram}}\kern-30pt\end{equation}                        
\begin{equation}\label{seq:monadKi}\begin{diagram} 
   \OO(-1,0)^{k}&\rTo_{\vect{A-y \\
                             B-x \\
                             D   \\
                             D_2A}}^{\A(K_\infty)}
    &\col{\OO(-1,1)^{k}\\ \oplus\\ \OO^{k+3}}
    &\rTo_{\vect{x-B  &A-y &   C    &0\\
                 -D_2 & 0  & \vect{0&-y} &1}}^{\B(K_\infty)}& \OO(0,1)^{k+1},
\end{diagram}\end{equation}                        
\begin{equation}\label{seq:monadK0i}\begin{gathered}
                                   \begin{diagram} 
   \OO(-1,0)^{k+j}
            &\rTo_{\vect{A-y         &  0\\
                         A'           & -y        \\
                         B-x          & -C_1e_+   \\
                         \vect{B'\\0} &s-x-C_1'e_+\\
                         D_2          &0          \\
                         0            &e_+}}^{\A(K_{0\infty}(0,1))}
    &\kern-12pt\col{\OO(-1,1)^{k+j}\kern-8pt\\ \kern-4pt\oplus \\ \OO^{k+j+2}}\kern-11pt
     &&&\rTo_{}^{\B(K_{0\infty}(0,1))\kern-10pt}&&
    &\kern-3pt \OO(0,1)^{k+j},\\
\end{diagram}\kern-15pt\\
\B(K_{0\infty}(0,1))=\vect{x-B           &C_1e_+      &A-y &0  &C_2 &AC_1 \\
                  \vect{-B'\\0} &x-s+C_1'e_+ & A' &-y &C_2'&A'C_1},
\end{gathered}\end{equation}                        
with $A, B, C, D, A', B'$ and $C'$ being matrices of respective size
$k\times k,k\times k, k\times 2, 2\times k, j\times k, 1\times k,
j\times 2$, ($C_i$, $C_i'$ being the $i$th column of $C_i$, $C_i'$, and
 $D_i$ the $i$th row of $D$), and satisfying the
monad equations
\begin{align}
  [A,B] + CD&=0,                                           \label{monadeqn1}\\
  \vect{B'\\0}A + sA' - A'B -C'D&=0,  \label{monadeqn2}\\
  -e_+ A'+ \vect{1&0} D &=0,                                    \label{monadeqn3}
\end{align}
and the genericity conditions
\begin{align}
\vect{A-y\\B-x\\D}
            &\text{ injective for all }x,y\in \C,\label{gencon1}\\
\vect{x-B&A-y& C}
            &\text{ surjective for all }x,y\in \C,\label{gencon2}\\
 \vect{x-B&A& C& 0\\ \vect{-B'\\0}& A'& C'& x-s\\
 0&0& \vect{1&0}& -e_+ }
            &\text{ surjective for all }x \in \C,\label{gencon3}
\\
N= \vect{\vect{A\\A'}&\vect{C_2\\ C'_2}&M\vect{C_2\\ C'_2}&\cdots&
           M^{j-1}\vect{C_2\\ C'_2}}&\text{ is an isomorphism},\label{gencon4}\\
\end{align}
where
\begin{equation}
M= \vect{ B& -C_1e_+   \\
                         \vect{B'\\0} &s-C_1'e_+},
\end{equation}
  modulo the action of $Gl(k,\C)$
\begin{equation}\label{action}
(A,B,C,D,A',B',C')\mapsto
        (gAg^{-1}, gBg^{-1}, gC, Dg^{-1},A'g^{-1}, B' g^{-1},C').
\end {equation}

The maps between the bundles are mediated by the following maps of monads
\begin{equation}\label{diag:mediated}
\begin{diagram}[tight,height=1.5em,width=6.7em]
\OO(-1,0)^{k+j}     &\rTo_{\A(K_0)}     &\OO(-1,1)^{k+j}\oplus\OO^{k+j+3}                     &\rTo_{\B(K_0)}       &&O(0,1)^{k+j+1}\\  \\ \\
\dTo_{\vect{1&0\\0&1}}&                    & \dTo_{\vect{1&0&0&0&0&0&0\\
                                                              0&1&0&0&0&0&0\\
                                                              0&0&1&0&0&0&-C_1\\
                                                              0&0&0&1&0&0&-C_1'\\
                                                              0&0&0&0&0&1&0\\
                                                              0&0&0&0&0&0&1}}                        &                     && \dTo^{\vect{1&0&-C_1\\0&1&-C_1'}}\\ \\ \\
\OO(-1,0)^{k+j}&\rTo^{\A(K_{0\infty}(0,1))\negthinspace\negthinspace\negthinspace\negthinspace}
                                        &\OO(-1,1)^{k+j}\oplus\OO^{k+j+2}                  &\rTo^{\B(K_{0\infty}(0,1))}       &&O(0,1)^{k+j}\\ \\ \\
\uTo_{\vect{A\\ A'}}&                   & \uTo_{\vect{A&0&0&0\\
                                                      A'&0&0&0\\
                                                      0&A&\vect{0&C_2}&0\\
                                                      0&A'&\vect{0&C_2'}&0\\
                                                      0&0&\vect{0&0\\1&0}&\vect{1\\0}}}    &                    & & \uTo^{\vect{A&C_2\\ A'&C_2'}}\\ \\ \\
\OO(-1,0)^{k}       &\rTo_{\A(K_\infty)}&\OO(-1,1)^{k}\oplus\OO^{k+3}                         &\rTo_{\B(K_\infty)}  &&O(0,1)^{k+1}\\ \\
\dTo_{\vect{1}}     &                   &\dTo_{\vect{1&0&0&0\\0&1&0&0\\0&0&1&0}}              &                     && \dTo^{\vect{1&0}}\\ \\
\OO(-1,0)^{k}       &\rTo^{\A(E)}       &\OO(-1,1)^{k}\oplus \OO^{k+2}                        &\rTo^{\B(E)}         &&O(0,1)^{k}\\ \\
\uTo_{\vect{1&0}}   &                   & \uTo_{\vect{1&0&0&0&0&0\\0&0&1&0&0&0\\0&0&0&0&1&0}} &                     && \uTo^{\vect{1&0&0}}\\ \\
\OO(-1,0)^{k+j}     &\rTo^{\A(K_0)}     &\OO(-1,1)^{k+j}\oplus\OO^{k+j+3}                   &\rTo^{\B(K_0)}       &&O(0,1)^{k+j+1}.
\end{diagram}\kern-20pt
\end{equation}

When $j=0$, the bundle $E$ and $K_\infty$ are the cohomology of the same monads
given by Equation (\ref{seq:monadE}) and (\ref{seq:monadKi}), with the matrices
$A,B,C,D$ satisfying Equations (\ref{monadeqn1}), (\ref{gencon1}), and (\ref{gencon2}), 
the matrix $A$ is invertible, and
the matrices for $K_0$ and $K_{0\infty}(0,1)$ are
\begin{align*}
\A(K_0)&=\vect{A-y\\ B-x\\ D\\ D_1A^{-1}},  &\B(K_0)&=\vect{x-B & A-y & C & 0\\
                                                         -D_1A^{-1}&0&\vect{1&0}&-y},\\
\A(K_{0\infty}(0,1))&=\vect{A-y\\ B-C_1D_1A^{-1}-x\\ \vect{D_2\\ D_1A^{-1}}}, 
           & \B(K_{0\infty}(0,1))&=\vect{x-B+C_1D_1A^{-1} & A-y & \vect{C_2&AC_1}}.
\end{align*} 
The maps between the bundle are mediated by a the $j=0$ version of Diagram (\ref{diag:mediated}).
\end{theorem}

The proof of this theorem is postponed to Section \ref{sec:proofMonads}.

\begin{remark}\label{remark:flag}
A flag of degree $j$ in $E|_{\PP^1\times\{0\}}$ is given by the projective equivalence
class of a
(pointwise) injective map of bundle
\[\OO(-j)\hookrightarrow E|_{\PP^1\times\{0\}}.\]
The bundle $E|_{\PP^1\times\{0\}}$ on $\PP^1$ splits
as a sum $\OO(n)\oplus\OO(-n)$ for some $n\in\N$.  The injection $\OO(-j)\hookrightarrow\OO(n)\oplus\OO(-n)$ is equivalent to a nowhere vanishing
section of $\OO(n+j)\oplus\OO(-n+j)$.
%
%
Thus the existence of a flag
$\OO(-j)\hookrightarrow E|_{\PP^1\times\{0\}}$ guarantees that
$E|_{\PP^1\times\{0\}}$ splits as $\OO(-j)\oplus \OO(j)$ when $j\leq 0$,
and as $\OO(n)\oplus \OO(-n)$ for some $0\leq n\leq j$ when $j\geq 0$.
Obviously, only $j\geq0$ matters for studying $E$, but it turns out the
result for $j\leq 0$ is useful for studying $K_0$.

The splitting of $E_0$ imposed by the existence of the
flag forces $\dim\ker(A)=n$.  Indeed, restrict  Monad
(\ref{seq:monadE})
to $\PP^1\times\{0\}$, and tensor with $\OO(-1)$ throughout to obtain $E|(-1)$
as the cohomology of the monad
\begin{diagram}
\OO(-2)^k&\rTo^\alpha &\OO(-2)^k\oplus\OO(-1)^{k+2} &\rTo^\beta &\OO(-1)^k.
\end{diagram}

One then has $h^1( \PP^1\times\{0\}, E(-1)) = {\rm dim\ coker}(A)$, hence
\begin{equation}\label{eqn:dimker}
    \dim\ker(A)=n\quad (\text{or $|j|$, if }j\leq0).\end{equation}
In particular, $A$ is invertible when $j=0$.
\end{remark}

Let us now consider the converse to Theorem \ref{thm:monadE}.
Given matrices satisfying the monad equations (\ref{monadeqn1}),
(\ref{monadeqn2}) and (\ref{monadeqn3}),  we
can construct cohomology sheaves $E$, $K_0$, $K_\infty$, $K_{0\infty}(0,1)$
with  maps $K_0\rightarrow E, K_\infty\rightarrow E, K_0\rightarrow
K_{0\infty}(0,1), K_\infty\rightarrow K_{0\infty}(0,1)$. The genericity
conditions (\ref{gencon1}), (\ref{gencon2}), and (\ref{gencon3}) ensures those
sheaves are bundles.
It is
routine, but lengthy, work to verify that those maps are
injective sheaf maps. Similarly, it is not hard to
verify that the cokernels of $K_0\rightarrow K_{0\infty}(0,1)$ and
$K_\infty\rightarrow E$ are $\OO_\infty$, as wanted. There remains the cokernels of the maps
$K_\infty\rightarrow K_{0\infty}(0,1)$ and
$K_0 \rightarrow E$.

 Suppose $f\in\ker(\B(E))$, with components $f_1,f_2,f_3$ in the
decomposition $\C^k\oplus\C^j\oplus\C^2$ used in Theorem \ref{thm:monadE}.
Away from $y=0$, $f$ is the image
of \[\vect{f_1\\ 0\\ f_2\\ \Bigl(-\vect{B'\\ 0}f_1+A'f_2+C'f_3\Bigr)/y
     \\ f_3\\ \vect{1&0}f_3/y}\in\ker(\B(K_0)),\]
hence the cokernel $Q_0$ of the map $K_0\to E$ is supported on
$y=0$.  To verify that
$Q_0=\OO(j,0)_0$, as in Equation (\ref{seq:defK0}), we verify
that $\pi_*Q_0=\OO(j)$, using the following lemma.

\begin{lemma}\label{lemma:direct-image-resolution}
Suppose that a bundle $F$ is defined by a monad
\begin{equation}\label{seq:monadtype}
\begin{diagram}
    \OO(-1,0)^{\ra}&\rTo_{\begin{bmatrix}A_1y+A_0\\B_1x+B_0\end{bmatrix}}
    &\col{\OO(-1,1)^{\rb}\\ \oplus\\ \OO^{\rc}}
    &\rTo_{\vect{C_1x+C_0&D_1y+D_0}}
     &\OO(0,1)^{\rd},\end{diagram}\end{equation}
and that $\pi_*F(0,-1)=0$. Then $\Ro F(0,-1)$ has a resolution
\begin{equation}\label{direct-image-resolution}
\begin{diagram} 0&\rTo&\OO(-1)^{\rb}&\rTo_{\vect{C_1x+C_0}} &\OO^{\rd}&\rTo&\Ro F(0,-1)&\rTo &0.\end{diagram}\end{equation}
\end{lemma}

\begin{proof}Use
K\"unneth Theorem for sheaves (\cite[Prop 9.2.4 p. 116]{BKkempf})
to see that
\[R^i\pi_*\OO(p,q)=\OO(p)\otimes H^i\bigl(\PP^1,\OO(q)\bigr).\]
Alternatively, use the projection formula \cite[Ex. 8.3, p.
253]{Hartshorne1977}, and   prove
$R^i\pi_*\OO(0,q)=\OO^{h^i(q)}$ using the exact sequence
$\OO(0,q)\to\OO(0,q+1)\to \OO|_{x=\infty}$.

Tensor Monad (\ref{seq:monadtype}) by $\OO(0,-1)$, and let
$K=\ker\vect{C_1x+C_0& D_1x+D_0}$. The short exact
sequence $0\to\OO(-1,-1)^{\ra}\to K\to F(0,-1) \to 0$ on $\PPP$
gives on $\PP^1$ the isomorphisms
\[R^i\pi_*K=R^i\pi_*F(0,-1), \text{ for }i\geq 0.\]

The short exact sequence $0\to K\to \OO(-1,0)^\rb\oplus
\OO(0,-1)^\rc\to \OO^\rd\to 0$ induces the exact sequence
\begin{diagram}
0&\rTo&\pi_*K&\rTo&\OO(-1)^{\ra}
&\rTo&\OO^{\rd}&\rTo&R^1\pi_*K&\rTo&0\end{diagram} on $\PP^1$, and
$R^i\pi_* K =0 \text{ for }i\geq 2$.
\theproofcomplete
\end{proof}

This lemma gives resolutions for $R^1\pi_*K_0(0,-1)$ and $R^1\pi_*E(0,-1)$,
with kernels the zero sheaves $\pi_*K_0(0,-1)$ and $\pi_*E(0,-1)$.     We then have a resolution diagram
\begin{diagram}
0&\rTo&\OO(-1)^{k\phantom{+j}}&\kern-10pt\rTo^{x-B}&\OO^{k}&\rTo&R^1\pi_*E(0,-1)&\rTo&0\\
&&\uTo&&\uTo&&\uTo\\
0&\rTo&\OO(-1)^{k+j}&\rTo^{\vect{x-B&0\\\vect{-B'\\0}&x-s\\0&-e_+}}&\OO^{k+j+1}&\rTo&R^1\pi_*K_0(0,-1)&\rTo&0\\
&&\uTo&&\uTo&&\uTo\\
0&\rTo&\OO(-1)^{j\phantom{+j}}&\kern-10pt\rTo^{\vect{x-s\\-e_+}}&\OO^{j+1}&\rTo&K&\rTo&0.
\end{diagram}
Diagram chasing gives us a kernel $K$, with its resolution; however, the last row is a fairly standard resolution of $\OO(j)$, and so $K=\OO(j)$.

It remains to verify that the cokernel $Q_\infty$ of the map
$K_\infty\rightarrow K_{0\infty}(0,1)$ is $\OO(-j,0)_0$. We first check that the map is surjective away from $y=0$.
Suppose $f\in\ker(\B(K_{0\infty}(0,1)))$, with components $f_1, \ldots, f_6$ in
the decomposition $\C^k\oplus\C^j\oplus\C^k\oplus\C^j\oplus\C\oplus\C$ used in
Theorem \ref{thm:monadE}.  Since
$\tiny \vect{A\\ A'}$ is injective, there exists $\vect{P& Q}$
such that $PA+QA'=1$.  Set
\[
\begin{aligned}
g_1&:=Pf_1+Qf_2,\\
g_3&:=f_6,\\
g_5&:=f_5,\\
g_4&:=\frac{g_5-D_2g_1}y,\\
g_2&:=Pf_3+Qf_4-(PC_2+QC_2')g_4,
\end{aligned}\quad\quad\quad\quad \text{ and }g:=\vect{g_1\\ \vdots\\ g_5}.\]
Away from $y=0$, we have $g$ is mapped to $f$, and
$g\in\ker(\B(K_\infty))$.  To verify this last statement, the only difficulty
is proving that
\[(x-B)g_1+(A-y)g_2+C_1g_3+C_2g_4=0.\]
The trick is to prove that
$A$ and $A'$ times the left-hand-side of the equation is $0$ and then use the injectivity of $\tiny\vect{A \\ A'}$.

Thus the sheaf $Q_\infty$ is supported over the
line $y=0$.

In the special case $j=0$, the components $f_2$ and $f_4$ are
automatically zero.
Since $A$ is invertible, we can
normalize in a neighborhood of $y=0$ to $f_1=0$.  Then $g_2=f_5/y$.
The only problem is when $f_5\neq0 $.  Hence $Q_\infty=\OO_0$, as desired.

Suppose then that $j\neq 0$.
Applying Lemma \ref{lemma:direct-image-resolution}, we have a resolution
\begin{equation}\label{res-diagram}
\begin{diagram}[height=2.6em]
\OO(-1)^{k+j}&\rTo^{\vect{x-B&C_1e_+\\ \vect{-B'\\0}&(x-s)+C'_1e_+}}
   &\OO^{k+j}&\rTo& R^1\pi_*(K_{0\infty})\\
\uTo^{\vect{A\\A'}}&&\uTo_{\vect{A&C_2\\A'&C'_2}}&&\uTo\\
\OO(-1)^{k}&\rTo^{\vect{x-B\\-D_2}}&\OO^{k+1}&\rTo& R^1\pi_*(K_\infty(0,-1)).
\end{diagram}\end{equation}

Note that  $R^1\pi_*(K_{0\infty})$ is supported on points, away from $\infty$. We want to build up
a resolution of the twists $R^1\pi_*(K_\infty(0,-1))(\ell)$, and then show that $R^1\pi_*(K_\infty(0,-1))(j)$ maps to $R^1\pi_*(K_{0\infty})$ with kernel $\OO$.  Suppose for the moment
that we have a sheaf $F$ with resolution
\begin{diagram}
0&\rTo&\OO(-1)^a&\rTo^{\vect{x+\alpha\\ \beta}}&\OO^{a+1}&\rTo &F&\rTo &0.
\end{diagram}
The injection $\OO\to\OO^{a+1}$ in the last coordinate induces a map
$s_1\colon\OO\to F$ which is non-zero at $\infty$.

Let $s(\ell)$ and $e_+(\ell)$ be $\ell\times\ell$ and $1\times \ell$ versions  of
the matrices defined by Equation (\ref{def:seplus}).  The sheaf $\OO(\ell)$ has
resolution
\begin{diagram}
0&\rTo&\OO(-1)^\ell&\rTo^{\vect{x-s(\ell)\\ -e_+(\ell)}}&\OO^{\ell+1}&\rTo &\OO(\ell)&\rTo &0,
\end{diagram}
and the injection $\OO\to\OO^{\ell+1}$ in the first coordinate induces a map
$s_2\colon\OO\to\OO(\ell)$.

The maps $s_1$ and $s_2$ induce naturally
$\tilde s_1\colon \OO(\ell)\to F(\ell)$ and $\tilde s_2\colon F\to F(\ell)$.
We then have a short exact sequence
\begin{diagram}
0&\rTo&\OO&\rTo^{\vect{s_1\\ -s_2}}&F\oplus\OO(\ell)&\rTo^{\vect{\tilde s_2& \tilde s_1}} &F(\ell)&\rTo &0.
\end{diagram}
Using the snake lemma on the resolutions
\begin{diagram}
0&\rTo&0&\rTo&\OO&\rTo&\OO&\rTo &0\\
&&\dTo&&\dTo&&\dTo\\
0&\rTo&\OO^{a+\ell}&\rTo&\OO^{a+1+\ell+1}&\rTo &F\oplus\OO(\ell)&\rTo &0,
\end{diagram}
we get the resolution
\begin{diagram}
0&\rTo&\OO^{a+\ell}&\rTo_{\vect{x+\alpha&0\\ \vect{\beta\\0}&x-s(\ell)\\0&-e_+(\ell)}}&\OO^{a+\ell+1}&\rTo &F(\ell)&\rTo &0.
\end{diagram}

Setting, as above,
\begin{equation}\label{eqn:defM}
M := \vect{B&-C_1e_+\\ \vect{B'\\0}&s-C'_1e_+},
\end{equation}
and
\begin{align}
N' &:=\vect{\vect{A\\ A'}&\vect{C_2\\ C'_2}&
               M\vect{C_2\\ C'_2}&\cdots&M^{j-2}\vect{C_2\\ C'_2}},\\
N\ &:= \vect{\vect{A\\ A'}&\vect{C_2\\ C'_2}&
               M\vect{C_2\\ C'_2}&\cdots&M^{j-1}\vect{C_2\\ C'_2}}
    \label{eqn:defN}
\end{align}
(so that $M, N$ are $(k+j)\times (k+j)$ matrices, and $N'$ is
$(k+j)\times (k+j-1)$), we build up a twist by $\OO(j-1)$ of
Diagram (\ref{res-diagram})
\begin{equation}\label{res-diagram-j}
\begin{diagram}
\OO(-1)^{k+j}&\rTo^{x-M}
                     &\OO^{k+j}&\rTo& R^1\pi_*(K_{0\infty})\\
\uTo^{N'}&&\uTo_{N}&&\uTo\\
\OO(-1)^{k+j-1}&\rTo_{\vect{x-B&0\\
   \vect{-D_2\\0}&x-s(j-1)\\ 0&-e_+(j-1)}}
            &\OO^{k+j}&\rTo& R^1\pi_*(K_\infty(0,-1))(j-1).
\end{diagram}\end{equation}
The kernel $\pi_*Q(j-1)$ of the map
$R^1\pi_*(K_\infty(0,-1))(j-1)\rightarrow R^1\pi_*(K_{0\infty})$ is the
desired $\OO(-1)$ if and only if the induced map $N$ on sections
is an isomorphism, hence the genericity condition we have  imposed on the
matrix $N$ in our theorem above.

While Theorem \ref{thm:monadE} gives us the matrices starting from the
bundles created from the knowledge of the flag and trivialisations, we just
proved that the matrices give us the bundles $E$, $K_0$, $K_\infty$,
and $K_{0\infty}(0,1)$, from which we can extract the flag and trivialisations.
We can then end this section with a theorem,
summarising the passage from a bundle and
flag to their associated monad.

\begin{theorem}\label{thm:BundleMatrices}
There is an equivalence between

1) Vector bundles $E$ of rank two on $\PP^1\times \PP^1$, with
   $c_1(E) = 0, c_2(E) =k$.
trivialized along $\PP^1\times\{\infty\}\cup\{\infty\}\times\PP^1$,
and with a based flag $\phi\colon\OO(-j)\hookrightarrow E$  of degree
$j$ along $\PP^1\times\{0\}$ (up to non-zero scalar multiple),
with the basing condition $\phi
(\infty) (\OO(-j)) = {\rm span}(0,1)$

2) Matrices $A,B$ ($k\times k$), $C$ ($k\times 2$), $D$ ($2\times k$),
$A'$ ($j\times k$), $B'$ ($1\times k$), $C'$ ($j\times 2$),
satisfying the monad equations (\ref{monadeqn1}), (\ref{monadeqn2}),
(\ref{monadeqn3})
and the genericity conditions
(\ref{gencon1}), (\ref{gencon2}), (\ref{gencon3}), (\ref{gencon4}) modulo the action of $Gl(k,\C)$ given by Equation (\ref{action}).
\end{theorem}

\section{Normal forms: the proof of Theorem \ref{thm:monadE}.}
\label{sec:proofMonads}
We now prove Theorem \ref{thm:monadE} of page~\pageref{thm:monadE}.  To do
so, we normalize the monads given by Theorem \ref{prop:monadEgen} so that
in the end they are defined only up to a $Gl(k,\C)$ action.
To simplify the notation given by Equation (\ref{eqn:notationV}), set
\[\begin{aligned}\VE{i}&:=V_i(E),\quad \quad \quad &\VK{i}&:=V_i(K_0),\\
\VKi{i}&:=V_{i}(K_\infty),\quad\text{ and }\quad
                         &\VKK{i}&:=V_{i}(K_{0\infty}(0,1)).
\end{aligned}\]

{}We first normalize the monads for $E$. From Lemma \ref{lemma:restraints}, we know that if
\[\A(E)=\vect{\alpha_1y+\alpha_0\\ \beta_1x+\beta_0}\quad
\text{ and } \quad \B(E)=\vect{\mu_1x+\mu_0 & \nu_1y+\nu_0},\] then
$\alpha_1$ is an isomorphism, $\beta_1$ is injective and $\nu_1$ is
surjective.  We also know that $\mu_1$ is injective, and since
$k_2=k=k_4$, it must be an isomorphism.  {}From Lemma \ref{lemma:restraints},
we also have that $V_3=\beta_1(V_1)\oplus \ker(\nu_1)$.

Given any basis of $\VE1$ and $\ker(\nu_1)$, we can pick bases of
$\VE2=\alpha_1(\VE1)$, $\beta_1(\VE1)$ and $\VE4=\mu_1\alpha_1(\VE1)$ so that
\begin{gather*} \alpha_1=-1_{k\times k},  \quad \mu_1=1_{k\times k},\\
\beta_1=\vect{-1_{k\times k}\\ 0_{r\times k}},\quad
\nu_1=\vect{-1_{k\times k}& 0_{k\times r}}.\end{gather*}
So we get
\[\A(E)=\vect{A-y\\ B-x\\ D},
  \text{ and } \B(E)=\vect{x+\mu_0 &\nu_{00}-y & C}.\]
Since $0=\B(E)\A(E)=(\mu_0 A+\nu_{00} B+CD)+(A-\nu_{00})x-(\mu_0+B)y$,
we must have $\mu_0=-B$, $\nu_{00}=A$, and $[A,B]+CD=0$.

Notice that we did not use all the freedom we were given by Remark
\ref{remark:symmetries} of page \pageref{remark:symmetries}:  the basis of
$\C^k=\VE1$ is still totally arbitrary.  However, the basis of
 $\C^2=\ker(\nu_1)$ is induced by the trivialization of $E$ along
$\PP^1\times\{\infty\}$. The  residual freedom is $Gl(k)$.

Let us continue with $K_0$, and normalize its monad. Some of the normalization is
inherited from that of $E$.
We know the inclusion $K_0\to E$ gives surjections $\Phi_i\colon
\VK{i}\to \VE{i}$. Consider the diagram
\begin{equation}\label{monadEK0}
\begin{diagram}[height=2.7em]
\col{\VK1\\ \otimes\\ \OO(-1,0)} &\rTo^{\vect{\alpha_1^0 y+\alpha_0^0\\
\beta_1^0x+\beta_0^0}} &\col{\VK2\otimes \OO(-1,1)\\ \oplus\\
\VK3\otimes\OO} &\rTo^{\vect{\mu_1^0x+\mu_0^0 & \nu_1^0y+\nu_0^0}}
&\col{\VK4\\ \otimes\\ \OO(0,1)}\\
\dOnto^{\Phi_1} && \dOnto_{\Phi_{23}=\Phi_2\oplus\Phi_3} && \dOnto^{\Phi_4}\\
\col{\VE1\\ \otimes\\ \OO(-1,0)} &\rTo^{\vect{A-y\\ B-x\\ D}} &\col{\VE2\otimes
\OO(-1,1)\\ \oplus \\ \VE3\otimes\OO} &\rTo_{\vect{x-B&A-y& C}}
&\col{\VE4\\ \otimes\\ \OO(0,1).}
\end{diagram}
\end{equation}

{}From the coefficient of $xy$ in the monad equations for $K_0$ and $E$,
we get the commutative diagram
\begin{diagram}
        &     &\VK2&      &\rInto^{\mu_1^0}    &    &\VK4\\
        &\ruTo^{\alpha_1^0}_{\cong}&\vLine_{\Phi_2}  &      &        &\ruOnto<{-\nu_1^0}\\
\VK1&     &\HonV   &\rInto^{\beta_1^0}  &\VK3&    &\dOnto_{\Phi_4}\\
        &     &\dOnto  &              &\dOnto^{\Phi_3}  \\
\dOnto^{\Phi_1}  &     &\VE2  &\hLine^{\mu_1}_{\cong}&\VonH   &\rTo&\VE4\\
        &\ruTo^{\alpha_1}_{\cong}&        &      &       &\ruOnto_{-\nu_1}\\
\VE1  &     &\rInto_{\beta_1}    &      &\VE3.
\end{diagram}

Let $\bar V_1$ be a copy of $\VE1$ in $\VK1$ so that
$\Phi_1|_{\bar V_1}$ is an isomorphism.  Set \[Z_i:=\ker(\Phi_i).\]
We have
\begin{align}\VK1&=\bar V_1\oplus Z_1,\label{splitV1K0} \\
             \VK2&=\alpha_1^0(\bar V_1)\oplus Z_2.\label{splitV2K0}
\end{align}

{}From Lemma \ref{lemma:restraints}, $\alpha_1^0$ is an isomorphism,
$\beta_1^0$
 and $\mu_1^0$ are injective, and $\nu_1^0$ is surjective. While
$\nu_1\beta_1$ is an
isomorphism, hence $\VE3=\Im(\beta_1)\oplus\ker(\nu_1)$,
at the level of $K_0$, we still have
$\ker(\nu_1^0)\cap\Im(\beta_1^0)=\{0\}$ but the direct sum doesn't
fill all of $\VK3$.  Note that $\beta_1^0(\ker\Phi_1)\subset
\ker(\Phi_3)$.  But more important is that
$\beta_1^0(\bar V_1)\cap \ker(\Phi_3)=\{0\}$ because
$\beta_1$ is injective.

Restrict the monads in  Diagram (\ref{monadEK0}) to $\PP^1\times\{\infty\}$.
{}From the display of those monads, we get information
about the various $\Phi_i$.  First, from the exact sequence
$0\to \VK1\otimes\OO(-1)\to \ker(\B(K_0|))\to K_0|\to 0$ and its equivalent
for $E$, we find in cohomology that the map
$H^0(\ker(\B(K_0|)))\to H^0(\ker(\B(E|)))$ is injective because it is
exactly the map $H^0(K_0|)\to H^0(E|)$.  But from the
sequence \[0\to \ker(\B(K_0|))\to \VK2\otimes \OO(-1)\oplus \VK3\otimes\OO
\to \VK4\otimes\OO\to 0\] and its equivalent for $E$, we get the
identifications
$H^0(\ker(\B(K_0|)))=\ker(\nu_1^0)$ and $H^0(\ker(\B(E|)))=\ker(\nu_1)$
compatible with the $\Phi_i$,
it must be that the restriction of $\Phi_3$ gives an isomorphism
$\ker(\nu_1^0)\to\ker(\nu_1)$.

Let $L\subset Z_3$ be a one-dimensional complement to $\beta_1^0(Z_1)$ in
$Z_3$.  Then
\begin{equation}\label{splitV3K0}
   \VK3=\beta_1^0(\bar V_1)\oplus \beta_1^0(Z_1)\oplus \ker(\nu_1^0)\oplus L.
\end{equation}

Note that again, $\mu_1^0(Z_2)\subset Z_4$.  Thus
\begin{equation}\label{splitV4K0}
   \VK4=\mu_1^0\alpha_1^0(\bar V_1)\oplus \mu_1^0(Z_2)\oplus \nu_1^0(L).
\end{equation}

The basis we have for the $\VE{i}$ can be lifted to induce basis of
$\bar V_1$,  $\alpha_1^0(\bar V_1)$, $\beta_1^0(\bar V_1)$,
$\mu_1^0\alpha_1^0(\bar V_1)$ and $\ker(\nu_1^0)$.  We can then write
\begin{gather*}
\Phi_1=\vect{1_{k\times k} & 0_{k\times j}}=\Phi_2\\
\Phi_3=\vect{1_{k\times k} & 0_{k\times j} &0_{k\times 2} &0_{k\times 1} \\
             0_{2\times k} & 0_{2\times j} &1_{2\times 2} &0_{2\times 1}}
\quad
\Phi_4=\vect{1_{k\times k} & 0_{k\times j} &0_{k\times 1}},
\end{gather*}
as in Diagram (\ref{diag:mediated}).

Given any basis for $Z_1$ and $L$, we can pick basis of $Z_2$,
$\beta_1^0(Z_1)$, and $\mu_1^0(Z_1)$ so that
\newcommand{\mm}{\phantom{-}}
\begin{gather*}\alpha_1^0=\vect{-1 &\mm0\\ \mm0&-1},\quad
  \beta_1^0= \vect{-1& \mm0\\ \mm0&-1\\ \mm0&\mm0\\ \mm0&\mm0},\quad
  \mu_1^0=\vect{1&0\\ 0&1\\ 0&0},\\
  \nu_1^0=\vect{-1&\mm0&0&\mm0\\ \mm0&-1&0&\mm0\\ \mm0&\mm0&0&-1}.\end{gather*}

The coefficients of $x$ and $y$ in the monad equation for $K_0$ and the
commutativity of diagram (\ref{monadEK0}) force
\newcommand{\aq}{{A'}}\newcommand{\aw}{{\alpha_{01}^0}}
\newcommand{\mq}{\mu_{01}^0}\newcommand{\mw}{\mu_{02}^0}
\newcommand{\ma}{\mu_{03}^0}\newcommand{\ms}{\mu_{04}^0}
\newcommand{\nq}{C'}\newcommand{\nw}{\nu_{01}^0}
\newcommand{\na}{\nu_{02}^0}\newcommand{\ns}{\nu_{03}^0}
\begin{gather*}
\alpha_0^0=\vect{A &0 \\ \aq & \aw},\quad
\mu_0^0 =\vect{-B&0\\ \mq&\mw\\ \ma&\ms},\\
\beta_0^0=\vect{B&0\\-\mq&-\mw\\ D&0\\-\ma&-\ms},
\quad\nu_0^0=\vect{A &0   &C  &0\\
                  \aq&\aw&\nq&\nw\\
                    0&0  &\na&\ns}.
\end{gather*}


Restrict $K_0$ and $E$ to $y=\epsilon\neq 0$, take duals and tensor by
$\OO(-1)$.  We have
\begin{diagram}
\VK4^*\otimes\OO(-1)&\rTo&\VK3^*\otimes\OO(-1)\oplus\VK2^*\otimes\OO&\rTo&\VK1^*\otimes\OO\\
\uTo                &    &                    \uTo                  &   &\uTo\\
\VE4^*\otimes\OO(-1)&\rTo&\VE3^*\otimes\OO(-1)\oplus\VE2^*\otimes\OO&\rTo&\VE1^*\otimes\OO.
\end{diagram}
The isomorphism $E^*|(-1)\to K_0^*|(-1)$ is mediated by the $\Phi_i^*$.  From
the display of the monads, we have $H^i(\ker\B^*(K_0|))=H^i(K_0^*|(-1))$,
and similarly for $E$.  Hence we have the commutative diagram of exact
sequences
\begin{diagram}
0\to&H^0(K_0^*|(-1))&\rTo&\VK2^*         &\rTo_{\vect{A^*-\epsilon&\aq^*\\
                                                        0&\aw^*-\epsilon}} &\VK1^*         &\rTo    &H^1(K_0^*|(-1))&\to  0      \\
    &  \uTo_{\cong}  &    &\uTo^{\Phi_2^*}&                                 &\uTo_{\Phi_1^*}&        &\uTo_{\cong}                   \\
0\to&H^0(E^*|(-1))  &\rTo&\VE2^*         &\rTo_{A^*-\epsilon}              &\VE1^*         &\rTo    &H^1(E^*|(-1))  &\to  0,
\end{diagram}
hence\[\ker\vect{A-\epsilon&0\\ \aq&\aw-\epsilon}\to\ker(A-\epsilon)\]
is an isomorphism through $\Phi_1$.  There can therefore be no kernel for
$\aw-\epsilon$ for all $\epsilon\neq0$.  Hence $\aw=0$.

Now, restricted to $y=0$, the sequence $0\to K_0\to E\to E_0/\Eof\to 0$ becomes,
for some $T$,
\[0\to T\to K_0|\to E|\to \OO(j) \to 0.\]  Since $c_1(K_0|)=0$,
we have $K_0|=\OO(m)\oplus\OO(-m)$.  The exact sequence forces $T=\OO(j)$.
Then the injection
$\OO(j)\to \OO(m)\oplus\OO(-m)$ forces $m=j$, as seen in Remark
\ref{remark:flag} on page \pageref{remark:flag}.

{}From the monad of $K_0|^*(-1)$,
\begin{diagram}\OO(-1)\otimes \VK4^*&\rTo
&\OO\otimes \VK2^*\oplus\OO(-1)\otimes \VK3^*&\rTo^{\A^*}&\OO\otimes \VK1^*,
\end{diagram}
we find that
\[\C^j=H^1(K_0|^*(-1)))=H^1(\ker(\A^*))=\coker \vect{A^*&\aq^*\\0&\aw^*}.\]
Thus it must be that
\[\dim\ker \vect{A&0\\ \aq&\aw}=j,\]
whence
\begin{equation}\label{AAinjective}
       \vect{A\\ \aq}\text{ is injective}.\end{equation}

Since $\A(K_0)$ is injective for all
$(x,y)=(x,0)$, it must be that
\[\vect{\mw+x\\ \ms}\colon Z_1\to\alpha_1^0(Z_1)\oplus L\]
is injective for all $x$.  Owing to Lemma (\ref{lemma:cyclicity})
below, we can then choose a basis of $Z_1$ and $\alpha_1^0(Z_1)\oplus L$
such that $\mw=-s$ and $\ms=-e_+$;
recall Equation (\ref{def:seplus}).

\begin{lemma}[Cyclicity] \label{lemma:cyclicity}
Suppose $\vect{T-\lambda\\ v}$ is an
injective map $\C^d\to \C^{d+1}$ for all $\lambda$.  Then
\[(vT^{d-1},\ldots,vT,v)\] is a basis for $\C^d$.
\end{lemma}

\begin{proof}
The result is invariant under conjugation $(T,v)\mapsto
(PTP^{-1},vP^{-1})$ and translation $(T,v)\mapsto (T-\lambda,v)$, hence we
only need to prove it for a matrix $T$ in Jordan normal form, with one
eigenvalue zero.  In fact, the injectivity hypothesis forces all the blocks
to have different eigenvalues.
We can then finish the proof by induction on the number of
Jordan blocks.
\end{proof}

Now that $\mu_{02}^0=-s$ and $\mu_{04}^0=-e_+$, we can perform
elementary column operations on $\beta_0^0$ to kill all but the first
line $B'$ of $-\mu_{01}^0$ and all of $-\mu_{03}^0$.  Such column
operations correspond to right multiplying by a matrix of the type
\[\vect{1&0\\ *& 1}\colon \bar V_1\oplus Z_1\to \bar V_1\oplus Z_1,\]
hence a repositioning of $\bar V_1$ in $\VK1$, while keeping $Z_1$ fixed.

Consider now the constant term of  $\B(K_0)\A(K_0)=0$.
Due to the splitting of $\VK1$ and $\VK4$, we find six equations, three
being Equations (\ref{monadeqn1}), (\ref{monadeqn2}), and (\ref{monadeqn3}),
one being tautologically $0$, and the remaining two being
\begin{align*}
 0=\nu_{02}^0e_+&\colon Z_1\to L\to \mu_0^0(Z_2),\\
 0=\nu_{04}^0e_+&\colon Z_1\to L\to \nu_1^0(L).
\end{align*}
Hence $\nu_{02}^0=0$ and $\nu_{04}^0=0$.  Once we consider that the
flag at $(\infty,0)$ lives in the second vector of the trivialization,
we have $\nu_{03}^0=\vect{1&0}$ in the appropriate basis.

We thus reduced the symmetries enough to establish the validity of
Equation (\ref{seq:monadK0}), and of the
fourth row of vertical maps in Diagram (\ref{diag:mediated}). Thus the monad for $K_0$ is as advertised, and the residual symmetry is $Gl(k)$, isomorphic to the symmetry of the monad for $E$.

Let us now continue with $K_\infty$,  with an obvious translate in the notation.
The problem is simplified as $\bar Z_1=\bar Z_2=\{0\}$, and $\dim
\bar Z_3=\dim \bar Z_4=1$.  We thus have, for some lift $F$ of $\ker(\nu_1)$,
\begin{align*}\VKi2&=\alpha_1^\infty(\VKi1),\\
              \VKi3&=\beta_1^\infty(\VKi1)\oplus F\oplus \bar Z_3,\\
              \VKi4&=\mu_1^\infty\alpha_1^\infty(\VKi1)\oplus \bar Z_4.
\end{align*}

Note that contrary to what happened for $K_0$, we cannot choose $F$ to
be $\ker(\nu_1^\infty)$, as it contains $\bar Z_3$. Indeed, the
exact sequence (\ref{seq:defKi})
for $K_\infty$, restricted to $y=\infty$, becomes
$0\to \OO\to K_\infty| \to E|\to\OO\to 0$, whence the map
$H^0(K_0|)\to H^0(E|)$ has a one-dimensional kernel.  Going through
the same analysis as before, where $H^0(K_0|)=\ker(\nu_1^\infty)$
and $H^0(E|)=\ker(\nu_1)$, we see that $\bar \Phi_3$ restricted to
$\ker(\nu_1^\infty)$ has a one-dimensional kernel, $\bar Z_3$ itself, as
claimed.

We lift the basis of the $\VE{i}$ to induce basis on all those pieces
of the $\VKi{i}$.  We can then write
\begin{gather*}\begin{aligned}
  &\bar\Phi_1=1_k & \bar\Phi_2&=1_k\\
  &\bar\Phi_3=\begin{bmatrix}1_k &0&0\\0&1_2&0\end{bmatrix}\quad
  &\bar\Phi_4&=\vect{1_k&0},
\end{aligned}\\
\alpha_1^\infty=-1_k\quad\quad
\beta_1^\infty=\vect{-1_k\\0\\0}\quad\quad
\mu_1^\infty=\vect{1_k\\0}\\
\nu_1^\infty=\vect{-1_k&0&0\\0&\vect{0&-1}&0}.
\end{gather*}

The commutativity of a diagram for $K_\infty$ analogous to Diagram
(\ref{monadEK0}) and the coefficients of $x$ and $y$ in the monad
equation $\B(K_\infty)\A(K_\infty)=0$ ensure
\begin{gather*}
\alpha_0^\infty=A,\quad \mu_0^\infty=\vect{-B\\-D_2},\quad
   \beta_0^\infty=\vect{B\\ D\\ \beta_{01}^\infty},\\
\nu_0^\infty=\vect{A&C&0\\0&\nu_{01}^\infty&\nu_{02}^\infty}.
\end{gather*}

Restricting the map $K_\infty\to E$ at $y=0$, where it is an isomorphism,
we have at the level of the cohomology of the monads that projection
on the first two factors must be an isomorphism
\[\ker\vect{A&C&0\\0&\nu_{01}^\infty&\nu_{02}^\infty}\to\ker\vect{A&C}.\]
For the projection to be an isomorphism, it must be that
$\nu_{02}^\infty\neq0$, and by choosing the basis of $\bar Z_3$ and $\bar Z_4$,
we can have $\nu_{02}^\infty=1$.  With column operations, we can kill
$\nu_{01}^\infty$, hence repositioning $F$ inside $\VKi3$.
The constant term of the monad equation then implies $\beta_{01}=D_2A$.

We thus established the validity of Equation (\ref{seq:monadKi}),
and the third row of  vertical maps in Diagram
(\ref{diag:mediated}).Thus the monad for $K_\infty$ is as announced, with the
same residual $Gl(k)$ symmetry.


Let us continue with $K_{0\infty}(0,1)$.
Notice that $K_{0\infty}(0,1)$ is to $K_0$ what $E$ is to
$K_\infty$.  Indeed, in a small neighborhood $U$ intersecting $y=\infty$,
\[\begin{aligned}
E(U)=&\OO(U)\oplus\OO(U), \text{ and }&K_\infty(U)=&\frac{\OO(U)}y\oplus\OO(U),&\text{ while}\\
K_{0\infty}(0,1)(U)=&y\OO(U)\oplus\OO(U), \text{ and }&K_0(U)=&\frac{y\OO(U)}y\oplus\OO(U).&
\end{aligned}\]
Also, $K_{0\infty}(0,1)$ is trivial on
$\PP^1\times\{\infty\}\cup\{\infty\}\times\PP^1$, and has a choice
of a flag in the section at $\infty$.

We can then use the monad of $K_0$ to extract the monad of $K_{0\infty}(0,1)$,
once however we normalize it correctly.  Staring at the monad given by
Equation (\ref{seq:monadKi}), we see we want an expression for $K_0$ of the
type
\[\A_2(K_0)=\vect{\tilde A-y\\ \tilde B -x\\ \tilde D\\ \tilde D_2\tilde A}, \quad
  \B_2(K_0)=\vect{x-\tilde B & \tilde A -y & \tilde C & 0\\
                                -\tilde D_2 & 0 & \vect{0& -y} & 1}.\]

To get an expression of this type, we set
\[P:=\vect{1&0&-C_1\\ 0&1&-C_1'\\ 0&0&1}, \quad
  Q:=\vect{ 1&0&0&0&0&0      &0\\
        0&1&0&0&0&0      &0\\
        0&0&1&0&0&C_1 &0\\
        0&0&0&1&0&C_1'&0\\
        0&0&0&0&0&0      &1\\
        0&0&0&0&1&0      &0\\
        0&0&0&0&0&1      &0},\]
and
\[\A_2(K_0):=Q^{-1}\A(K_0),\quad \B_2(K_0):=P\B(K_0)Q.\]

Now deleting the last row and column of $\B_2(K_0)$ and the last row of
$\A_2(K_0)$, we obtain the monad of Equation (\ref{seq:monadK0i}), and
establish the validity of the first row of vertical maps in Diagram
(\ref{diag:mediated}).

The only part of Theorem \ref{thm:monadE} that remains to be proved is the
validity of the second row of vertical maps in Diagram (\ref{diag:mediated}).  Notice
the map from sections of $K_\infty$ to sections $K_{0\infty}(0,1)$ is
multiplication by $y$, as  $K_{0\infty}(0,1)$ equals $K_{0\infty}(\PP^1\times\{\infty\})$,
not $K_{0\infty}(\PP^1\times\{0\})$.

On the dense set
$\{(x,y)\in\PPP\mid y\neq0\text{ and }(A-y) \text{ is invertible}\}$,
we can trivialize
the bundles by sending $\gamma\in\C^2$ to
\begin{gather*}
\vect{0\\-(A-y)^{-1}C\gamma\\ \gamma},  \quad
 \vect{0\\-(A-y)^{-1}C\gamma\\ \gamma\\ y\gamma_2},\quad
\vect{0\\0\\ -(A-y)^{-1}C\gamma\\ -y^{-1}A'(A-y)^{-1}C\gamma+y^{-1}C'\gamma\\ \gamma\\ y^{-1}\gamma_1},\\
\vect{0\\0\\-(A-y)^{-1}(y^{-1}AC_1\gamma_1+C_2\gamma_2)\\
          y^{-1}(y^{-1}A'C_1\gamma_1+C_2'\gamma_2)
                 -y^{-1}A'(A-y)^{-1}(AC_1\gamma_1y^{-1}+C_2\gamma_2)\\
          \gamma_2\\ y^{-1}\gamma_1}
\end{gather*}
respectively for $E$, $K_\infty$, $K_0$, and $K_{0\infty}(0,1)$.  This choice of trivialization
is preserved by the various $\Phi_{23}$ of Diagram (\ref{diag:mediated}) whose
validity we already established. For the proposed $\Phi_{23}$ of
$K_\infty\to K_{0\infty}(0,1)$, we have, using an obvious notation,
\[\Phi_{23}(\gamma^{K_\infty})=y\gamma^{K_{0\infty}(0,1)}.\]
Since the candidate $\Phi_{23}$ is globally defined and agrees with the actual
$\Phi_{23}$ on a dense subset of $\PPP$,  its validity is established.
The commutativity of the diagram forces $\Phi_1$ and $\Phi_4$ to be as claimed
in the second row of Diagram (\ref{diag:mediated}).

The genericity conditions are simply those for monads implied by Buchdahl's
Beilinson's theorem, along with the constraint on the matrix $N$ which was proven in the previous section.

The proof of Theorem \ref{thm:monadE} is now complete.


\section{\texorpdfstring{From monads to sequences of sheaves on $\PP^1$, and back again.}{From monads to sequences of sheaves on PxP, and back again.}}
\label{sec:seqsheaves}
We show in this section how  Diagram (\ref{Nahm-cplx-diagram}),
encoding the bundle $E$ and the flag, with the trivialization
over $y=\infty$, leads one to a Nahm complex, and inversely how the
Nahm complex encodes the diagram. We already have that the diagram
gives the monads of Theorem \ref{thm:monadE}
and morphisms between them, and conversely,
monads of our normalized form give back the diagram of bundles. It
thus suffices to show how our monad matrices encode, and are
encoded by, a Nahm complex.

The intermediary step we introduce are the exact sequences
(\ref{direct-image-sequences}) of Lemma \ref{lemma-on-direct-images}.
More specifically, set
\begin{equation}\label{sheavesonP1}\begin{aligned}
P_0 &:=R^1\pi_*(K_{0}(0,-1)),\\
P_\infty&:=R^1\pi_*(K_{\infty}(0,-1)),\\
Q_{\infty0}&:=R^1\pi_*(E(0,-1)), \\
Q_{0\infty}&:=R^1\pi_*(K_{0\infty}).
\end{aligned}\end{equation}
Those exact sequences can now be read
\begin{equation}\label{sheaf-sequences}
\begin{matrix}
0&\rightarrow& \OO(j)&\rightarrow& P_0&\rightarrow& Q_{\infty0}&\rightarrow&0,\cr
0&\rightarrow& \OO&\rightarrow& P_\infty&\rightarrow& Q_{\infty0}&\rightarrow&0,\cr
0&\rightarrow& \OO&\rightarrow& P_0&\rightarrow& Q_{0\infty}&\rightarrow&0,\cr
0&\rightarrow& \OO(-j)&\rightarrow& P_\infty&\rightarrow& Q_{0\infty}&\rightarrow&0,\cr
\end{matrix}
\end{equation}
with $Q_{0\infty}, Q_{\infty0}$ torsion sheaves of length $k+j, k$
respectively, supported away from infinity, while $P_0, P_\infty$
are trivialized over infinity on the line.

We note that the sheaves come with resolutions, as given by Lemma
\ref{lemma:direct-image-resolution}. From Theorem \ref{thm:monadE},
in the case $j\geq 1$, we get a diagram of resolutions
\begin{equation}\label{resolution-diagram}
\begin{diagram}
 \OO(-1)^{k+j}&\rTo
                  &\OO^{k+j+1}&\rTo& P_0&\rTo &0\\
\dTo^{\vect{1&0\\0&1}}&&\dTo_{\vect{1&0&-C_1\\0&1&-C'_1}}&&\dTo&\\
\OO(-1)^{k+j}&\rTo^{\vect{x-B&C_1e_+\\ \vect{-B'\\0}&(x-s)+C'_1e_+}}&\OO^{k+j}&\rTo& Q_{0\infty}&\rTo &0\\
\uTo^{\vect{A\\A'}}&&\uTo_{\vect{A&C_2\\A'&C'_2}}&&\uTo&\\
\OO(-1)^{k}&\rTo^{\vect{x-B\\-D_2}}&\OO^{k+1}&\rTo& P_\infty&\rTo &0\\
\dTo^{\vect{1}}&&\dTo_{\vect{1&0}}&&\dTo&\\
\OO(-1)^{k}&\rTo^{\vect{x-B}}&\OO^{k}&\rTo& Q_{\infty0}&\rTo &0\\
\uTo^{\vect{1&0}}&&\uTo_{\vect{1&0&0}}&&\uTo&\\
 \OO(-1)^{k+j}&\rTo^{\vect{x-B&0\\ \vect{-B'\\0}&x-s\\0&-e_+}}&\OO^{k+j+1}&\rTo& P_0&\rTo &0.\\
\end{diagram}\end{equation}

The next step is to prove that these resolutions always exist in this form, given the sheaves fitting into  Sequences  (\ref{sheaf-sequences}).
\begin{lemma}\label{lemma:fit}
Let $j\geq 1$ and let $P_0$, $P_\infty$, $Q_{0\infty}$, $Q_{\infty0}$ be sheaves over $\PP^1$
fitting into exact sequences (\ref{sheaf-sequences}), with $Q_{0\infty}$,
$Q_{\infty0}$ torsion of length $k+j, k$ respectively, supported away
from infinity in $\PP^1$. Then one has a commuting diagram of resolutions
as Diagram (\ref{resolution-diagram}).  Furthermore, the matrix $N$ defined by
Equation (\ref{eqn:defN}) using the matrices of the diagram is an isomorphism.
\end{lemma}

\begin{proof}
If one takes the resolution
\begin{diagram}
0&\rTo&\OO(-1,-1)&\rTo^{x-y}&\OO&\rTo&\OO_\Delta&\rTo&0,
\end{diagram}
of the diagonal in $\PP^1\times\PP^1$, lifts a  sheaf $F$ from $\PP^1$, tensors
 it with the resolution, and pushes down to the other factor, one obtains a
resolution
\begin{diagram}
0&\rTo&\col{H^0(\PP^1,F(-1))\\ \otimes\\ \OO(-1)}
   &\rTo^{\alpha_0(F)+\alpha_1(F)x}
               &\col{H^0(\PP^1,F)\\ \otimes\\ \OO}&\rTo&F&\rTo&0,\\
 &&&[f]\otimes g \mapsto [f]\otimes xg - [yf]\otimes g
\end{diagram}
as long as $H^1(\PP^1,F(-1))$ vanishes, which is the case for our
sheaves. Taking the sheaves and the maps $P_0\rightarrow
Q_{0\infty}\leftarrow P_\infty \rightarrow Q_{0\infty}\leftarrow
P_0$ and applying this process to them, we obtain a diagram akin to
Diagram (\ref{resolution-diagram}), with the right sheaves. One must show that the maps can be normalized as
advertised.

We first note that as $Q_{\infty0}$ is supported away from infinity,
one can identify \[H^0(\PP^1, Q_{\infty0}(-1)) = H^0(\PP^1,
Q_{\infty0})\] and normalize
$x\alpha_1(Q_{\infty0})+\alpha_0(Q_{\infty0})$ to $x- B$, for some
$B$.  The fact that $P_0,P_\infty$ are line bundles at infinity
allow us to filter their sections by order of vanishing at infinity.
We can split $H^0(\PP^1,P_0(-1))$, and $H^0(\PP^1,P_0)$ as  sums
\begin{equation*}
       H^0(\PP^1,Q_{\infty0}) \oplus H^0(\PP^1,\OO(j-1)),\text{ and }
       H^0(\PP^1,Q_{\infty0}) \oplus H^0(\PP^1,\OO(j)),
\end{equation*}
with the second summands being the kernels of projection to $Q_{\infty0}$, and
the first identifying $H^0(\PP^1,Q_{\infty0})$ with the subspace of sections
of $P_0(-1)$, $P_0$ vanishing at least to order $j,j+1$ at infinity. The
spaces $H^0(\PP^1,\OO(\ell)) $ have natural bases of sections
$1,\ldots,y^{\ell}$ in terms of which the resolution naturally
gets expressed in terms of the shift matrix $s$.
Finally, for a class $[f]\in H^0(\PP^1,Q_{\infty0})\subset H^0(\PP^1,P_0(-1))$,
we have
$[yf]\not\in \Im \bigl(H^0(\PP^1,\OO(j-1))\rTo^{y}H^0(\PP^1,\OO(j))\bigr)$,
hence the need for the matrix $B'$.
The last two lines of the diagram, and the maps between them, are then as
advertised.

Similarly for $P_\infty$, we can write
$H^0(\PP^1,P_\infty(-1))$, and $H^0(\PP^1,P_\infty)$ as
\begin{equation*}
H^0(\PP^1,Q_{\infty0}),\text{ and } H^0(\PP^1,Q_{\infty0}) \oplus H^0(\PP^1,\OO),
\end{equation*}
and show that the maps of the third row and between the third and fourth row
of the diagram is of the form given for some row vector $D_2$.

We then have isomorphisms $H^0(\PP^1,P_0(-1))\simeq H^0(\PP^1,Q_{0\infty}(-1))$
 showing us that we can take the map between the first elements of the first
and second rows to be the identity. The isomorphism
$H^0(\PP^1,Q_{0\infty}(-1))\simeq H^0(\PP^1,Q_{0\infty})$ shows us that we
can normalize $x\alpha_1(Q_{0\infty})+\alpha_0(Q_{0\infty})$ to $x-\tilde B$, for
some $\tilde B$.
The commutativity of the diagrams then tells us that $\tilde B$ is of
the form given, and that the remaining maps are also of the form given, for
suitable $C_1, C'_1,C_2, C'_2, A, A'$.

Use the monad equation (\ref{monadeqn1}) to set $D_1=e_+A'$.
The commutativity of Diagram (\ref{resolution-diagram})
implies the monad equations (\ref{monadeqn1}) and (\ref{monadeqn2}).

Finally, the fact that $N$ is an isomorphism follows from the fact that the map
$P_\infty(j)\rightarrow Q_{0\infty}$ must induce an isomorphism on sections, as
in Diagram (\ref{res-diagram-j}).
\theproofcomplete
\end{proof}

The genericity conditions on the matrices are equivalent to some
non-de\-ge\-ne\-ra\-cy conditions on the sheaves of Equation (\ref{sheavesonP1}).
We first note that one of the genericity conditions on the matrices is automatic
if they come from our sheaves. Indeed, if the condition (\ref{gencon3}) is not
satisfied (here $j>0$), we have that for some $x$,
\begin{equation}\label{surjective2}
\Im\vect{x-B&A& C& 0\\ \vect{-B'\\0}
   & A'& C'& x-s\\ 0&0& \vect{1&0}& -e_+ }\subset V
\end{equation}
with $V$ a proper codimension one subspace of $\C^{k+j+1}$.  Then
\[\begin{aligned}
\Im\vect{x-B&0\\ \vect{-B'\\0}&x-s\\ 0 &-e_+}
   &\subset V,
           \quad\Im\vect{A&C_2\\ A'&C_2'}\subset (V\cap\C^{k+j}),\\
\Im\vect{x-B& C_1e_+\\ -\vect{B'\\ 0}& x-s+C_1'e_+}
  &=\Im\Biggl(\Bigl(1-\vect{0&0&C_1\\ 0&0&C_1'\\0&0&1}\Bigr)
     \vect{x-B&0\\ \vect{-B'\\0}&x-s\\ 0 &-e_+}\Biggr)\\ &\subset V\cap \C^{k+j}.
\end{aligned}\]
Hence  by replacing
the spaces in the second column by
$V, V\cap\C^{k+j},\C^{k+1}, \C^k$, and $V$, we can ``reduce'' the diagram at $x$:
there are subsheaves $\tilde P_0,\tilde Q_{0\infty}$ which,
together with $P_\infty$ and $Q_{\infty0}$, fit in a variation of
Sequences (\ref{sheaf-sequences}) for which $j$ is replaced by $j-1$, and
giving as quotients of $P_0,Q_{0\infty},P_\infty,Q_{\infty0}$ the sheaves
$\C_x,\C_x,0,0$. In particular, the map $P_\infty\rightarrow Q_{0\infty}$ is not
surjective and we have left our class of sheaves.

The remaining conditions on our matrices do both translate into irreducibility for
our diagram. Let us say that the ``complex'' of sheaves is \emph{reducible} at $x$
if either
\begin{itemize}
\item {\it Case 1.}  There are skyscraper subsheaves $\C_x,\C_x,\C_x,\C_x$
of $P_0,Q_{0\infty},P_\infty,Q_{\infty0}$, localized at $x$, mapping to each other
by Sequences (\ref{sheaf-sequences}).

\item {\it Case 2.}  There are subsheaves $\tilde P_0,\tilde Q_{0\infty},\tilde P_\infty,\tilde Q_{\infty0}$ of $P_0,Q_{0\infty},P_\infty,Q_{\infty0}$, fitting in
Sequences (\ref{sheaf-sequences}), and giving as quotients the
sheaves $\C_x,\C_x,\C_x,\C_x$.
\end{itemize}

Translated to the world of resolutions, we can say that Diagram
(\ref{resolution-diagram}) is {\it reducible} at $x$ if either

\begin{itemize}
\item {\it Case 1.} There are, restricting at  $x$ so that
we are dealing with vector
spaces, one dimensional subspaces $V_1, V_2, V_3, V_4, V_5=V_1$ of the spaces in the
first column, (the subscript indicates the row) that lie in the kernel of the
maps from the first column to the second and that are mapped to each other
under the vertical maps; the
spaces in the first columns can then be replaced by quotients, giving a
``smaller'' diagram;

\item {\it Case 2.} There are, restricting at $x$, codimension one subspaces
$V_1, V_2,V_3,V_4$, $V_5=V_1$
of the spaces in the second column, containing the images of the maps between
the first and second column, and  mapped to each other under the vertical
maps, so that the diagram  can then be ``reduced'' to a smaller one.
\end{itemize}

We remark that it suffices to take dimension one or codimension one subspaces; other cases are reducible to this one, as can easily be checked.

We can then see that the   genericity conditions
(\ref{gencon1}), (\ref{gencon2}) on the matrices are equivalent
to Diagram (\ref{resolution-diagram}) being irreducible at  all $x$.

\begin{itemize}
\item {\it Case 1.} Suppose  there exists a one-dimensional subspace
$L\subset \C^k$ such that
\begin{equation}\label{injective}
L \subset \ker\vect{ A-y \\ B-x \\ D}
\end{equation}
for some $x$ and $y$. The monad equations then imply
\begin{equation*}\begin{aligned}
(B-x)(L)&= 0,           & \quad \Bigl(\vect{B'\\0}A + (s-x)A' \Bigr)(L) &= 0,\\
D_2(L)&= 0,             &  A(L)&\subset  L,\\
D_1(L)= e_+A'(L) &= 0.    &      &
\end{aligned}\end{equation*}

Hence there are subspaces
${\tiny\vect{A\\ A'}}L$, ${\tiny\vect{A\\ A'}}L$, $L$, $L$,
${\tiny\vect{A\\ A'}}L$
of the kernels in all the exact sequences in Diagram
(\ref{resolution-diagram}), which are mapped to each other under the vertical
maps. We can reduce the diagram at $x$.

\item {\it Case 2.} Suppose that for some $x$ and $y$,
\begin{equation}\label{surjective1}
\Im\vect{x-B&A-y& C}\subset V,
\end{equation}
with $V$ a proper codimension one subspace of $\C^{k}$, hence
\begin{equation*}\begin{aligned}
\Im(x-B)&\subset V,    &\quad \Im(C_1)&\subset V,\\
A(V)&\subset V,      &      \Im(C_2)&\subset V.
\end{aligned}\end{equation*}
Hence  by replacing
the spaces in the second column by
$V\oplus \C^{j+1}, V\oplus \C^j$,  $V\oplus \C, V$, and $V\oplus \C^{j+1}$,
we can again reduce the diagram at $x$
\end{itemize}

Conversely, if the diagram  is reducible at $x$, with a
common kernel through the diagram, one find an $L$ as in Equation
(\ref{injective}).

If the diagram is reducible at $x$, not with a common
kernel, but with a common one-dimensional cokernel, we have at $x$ codimension one subspaces $V_1,V_2,V_3,V_4$ of the spaces in
the second column, mapped to each other under the vertical maps and containing the images of the horizontal maps.  Because $\codim V_4=1$, there must be a
line $W\subset \C^k\subset \C^{k+1}$ so that $V_3\oplus W=\C^{k+1}$ and
$V_4\oplus W=\C^k$.  Since
$\vect{1&0&0}V_1\subset V_4$, we must have $V_1=V_4\oplus \C^{j+1}$.
To have $\codim V_2=1$, we must then have $V_2=V_4\oplus \C^j$.   Since
\[\vect{ A&C_2\\ A' & C_2'}(V_3)\subset V_2,\] we have
$\vect{ A & C_2}(V_3)\subset V_4$.
In the decomposition $V_4\oplus W\to V_3\oplus W$ we can write
\[\vect{ A & C_2}=\vect{ * & * \\ 0 & y}\]
for some $y\in\C$.   Then $\Im\vect{ A-y & C_2}\subset V_4$.  Since
$V_2$ contains the image of the second horizontal map, we must
have as well $\Im\vect{x-B&C_1}\subset V_4$.  But then
\[\Im\vect{A-y&x-B&C}\subset V_4\neq C^k,\]
hence the non-degeneracy condition (\ref{gencon2}) is not satisfied.

\begin{theorem}\label{thm:MatricesSequences}
There is an equivalence between

1) Matrices $A,B$ ($k\times k$), $C$ ($k\times 2$), $D$ ($2\times k$), $A'$ ($j\times k$),
$B'$ ($1\times k$), $C'$ ($j\times 2$),  satisfying the monad equations
(\ref{monadeqn1}), (\ref{monadeqn2}), (\ref{monadeqn3}),
and the genericity conditions
(\ref{gencon1}), (\ref{gencon2}), (\ref{gencon3}), (\ref{gencon4}) of page \pageref{gencon3},
modulo the action of $Gl(k,\C)$ given by Equation (\ref{action});

2) Exact sequences of sheaves
\begin{equation}\label{sheaf-sequences2}
\begin{matrix}
0&\rightarrow& \OO(j)&\rightarrow& P_0&\rightarrow& Q_{\infty0}&\rightarrow&0\cr
0&\rightarrow& \OO&\rightarrow& P_\infty&\rightarrow& Q_{\infty0}&\rightarrow&0\cr
0&\rightarrow& \OO&\rightarrow& P_0&\rightarrow& Q_{0\infty}&\rightarrow&0\cr
0&\rightarrow& \OO(-j)&\rightarrow& P_\infty&\rightarrow& Q_{0\infty}&\rightarrow&0
\end{matrix}
\end{equation}
on $\PP^1$, with $Q_{0\infty}, Q_{\infty0}$ torsion sheaves of
length $k+j, k$ respectively, supported away from infinity, and with
$P_0, P_\infty$ trivialized over infinity on the line, and irreducible.
\end{theorem}

\section{To Nahm complexes, and back again.}\label{sec:Nahmcomplex}
We next show that the sheaves fitting in the exact sequences
(\ref{sheaf-sequences2}) define,
and are defined by, a Nahm complex on the circle. We define these for the
integers $k>0$, $j\ge 0$,

If $j>0$, the Nahm complexes that we  consider over the circle are defined
by
\begin{itemize}
\item A bundle $V_{\infty0}$ of rank $k$ over the interval
$[\theta_0, \theta_\infty]$, equipped with a smooth connection
$\alpha_{\infty0}$, and a
covariant constant smooth section $\beta_{\infty0}$ of $End(V_{\infty0})$;
\item A bundle $V_{0\infty}$ of rank $k+j$ over the interval
$[\theta_\infty, 2\pi+ \theta_0]$, equipped with an smooth connection
$\alpha_{0\infty}$ on the interior, analytic near the boundary points, and a covariant constant section
$\beta_{0\infty}$ of $End(V_{0\infty})$ smooth on the interior, analytic near the boundary points;
\item At the boundary point $\theta_0$, an injection
$i_0\colon V_{\infty0}\rightarrow V_{0\infty}$ and a surjection
$\pi_0\colon V_{0\infty}\rightarrow V_{\infty0}$, such that
$\pi_0i_0 = Id$, so that one can decompose $V_{0\infty}$ as $
\ker(\pi_0)\oplus \Im(i_0)$. One asks that there be an extension of
this decomposition to a trivialization on the interior of the
interval such that one can write the connection $\alpha_{0\infty}$
and the endomorphism $\beta_{0\infty}$ in block form as
\begin{equation*}
\alpha_{0\infty}(t)=\vect{U(t)                & t^{\frac{j-1} {2}}W(t)\\
                         t^{\frac{j-1} {2}}V(t)&X(t)                 },\quad
\beta_{0\infty}(t)=\vect{P(t)                & t^{\frac{j-1} {2}}Q(t)\\
                        t^{\frac{j-1} {2}}R(t)&S(t)}
\end{equation*}
where $t$ is a local parameter with the point $\theta_0$ corresponding to
$t=0$. The top left blocks are $k\times k$, the bottom right block is
$j\times j$; $U,W,V,P,Q,R$ are analytic at $t=0$, and $X,S$ are meromorphic
with simple poles at $t=0$, and residues conjugate to
\begin{gather}
X_{-1}= diag({\frac{-(j-1)}{4}},{\frac{2-(j-1)}{4}},\dots,{\frac{(j-1)}{4}}),\\
S_{-1}= \vect{0&0&0&\dots&0&0\\ 1&0&0&\dots&0&0\\0&1&0&\dots&0&0\\
\dots&\dots&\dots&&\dots&\dots\\
0&0&0&\dots&1&0}.
\end{gather}
Furthermore,
\begin{equation*}
U(0) = \alpha_{\infty0}(0)\quad P(0) = \beta_{\infty0}
\end{equation*}
\item At the boundary point $\theta_\infty$, the boundary conditions are the
same as at $\theta_0$.
\item At both boundary points, some extra data, consisting of a
trivialization (choice of vectors $v_0, v_\infty$) of the $\frac{-(j-1)}{4}$
eigenspace of $X_{-1}$.
\end{itemize}

For $j= 0$, the constraints are simpler:
 the Nahm complexes over the circle that we consider are defined
by
\begin{itemize}
\item A bundle $V_{\infty0}$ of rank $k$ over the interval
$[\theta_0, \theta_\infty]$, equipped with a smooth connection
$\alpha_{\infty0}$, and a
covariant constant smooth section $\beta_{\infty0}$ of $End(V_{\infty0})$;
\item A bundle $V_{0\infty}$ of rank $k$ over the interval
$[\theta_\infty, 2\pi+ \theta_0]$, equipped with an smooth connection
$\alpha_{0\infty}$ and a covariant constant smooth section
$\beta_{0\infty}$ of $End(V_{0\infty})$
\item At the boundary point $\theta_0$, isomorphisms
$i_0\colon V_{\infty0}\rightarrow V_{0\infty}$, $\pi_0= i_0^{-1}$
with the gluing condition that
$\beta_{\infty0}-\pi_0\beta_{0\infty}i_0$ has rank one.
\item At the boundary point $\theta_\infty$, isomorphisms
$i_\infty\colon V_{\infty0}\rightarrow V_{0\infty}$, $\pi_\infty=
i_\infty^{-1}$ with the gluing condition that
$\beta_{\infty0}-\pi_\infty\beta_{0\infty}i_\infty$ has rank one.
\item At both boundary points, extra data consisting of decompositions $v_0 = (u_0, w_0)$, $v_\infty = (u_\infty,w_\infty)$ of the rank one boundary difference matrices  $\beta_{\infty0}-\pi_0\beta_{0\infty}i_0$, $\beta_{\infty0}-\pi_\infty\beta_{0\infty}i_\infty$ into products of a column and a row vector:
\begin{equation}
\beta_{\infty0}-\pi_0\beta_{0\infty}i_0= u_0\cdot w_0,\quad \beta_{\infty0}-\pi_\infty\beta_{0\infty}i_\infty = u_\infty\cdot w_\infty\end{equation}
\end{itemize}

There is a group $\mathcal G$ of gauge transformations which acts on the
complex and that can be used to normalize the complex as in the lemma below.
This group is constructed as follow:
one takes smooth $g_{\infty0}(z)\in Aut(V_{\infty0})$ on
$[\theta_0, \theta_\infty]$, $g_{0\infty}(z)\in Aut(V_{0\infty})$ on
$[\theta_\infty, 2\pi+ \theta_0]$ with, on the ``large'' side of the boundary
points, in the trivialisations used above, the constraint that
$g_{0\infty}(z)$ be analytic, with a decomposition
\begin{equation*}
g_{0\infty}(t)= \vect{K(t)& t^{\frac{j+1} {2}}L(t)\\
                     t^{\frac{j+1} {2}}M(t)&N(t)}
\end{equation*}
with $K, L, M, N$ analytic at $t=0$, and $K(0) = g_{\infty0}(0)$.
The group $\mathcal G$ acts as
\[g\cdot(\alpha,\beta) =
    (g\alpha g^{-1}-\frac12\dot{g}g^{-1},g\beta g^{-1}).\]

\begin{lemma}[Prop 1.15 of \cite{hurtubiseClassification}]
\label{lemma:normalform}
\begin{itemize}
\item Away from the boundary points, or even at the boundary points, if one is
on $V_{\infty0}$, one can gauge to $\alpha = 0, \beta =$ constant.

\item At the boundary, over $V_{0\infty}$, one can gauge transform to the
block form
\begin{align}
\alpha_{0\infty} &= \frac{1}{t}\vect{0&0\\0&diag({\frac{-(j-1)}{4}},{\frac{2-(j-1)}{4}},\dots,{\frac{(j-1)}{4}})},\\
\beta_{0\infty}&= \vect{ P_0&t^{\frac{j-1}{2}}q_0e_+\\
t^{\frac{j-1}{2}}\vect{r_0\\0}&-t^{-1}s + t^{j-1} \tilde s_0e_+}
\end{align}

Here $P_0$ is $k\times k$, $r_0$ is $1\times k$, $q_0 $ is $k\times 1$,
$\tilde s_0$ is $j\times 1$ and $P_0,r_0, q_0$ are constant in $t$, and,
setting $(\tilde s_0)_i= t^{i-1}(s_0)_i$, then $s_0$ is also constant
in $t$.

\item Using the  gauge transformation
\[G(t):=diag (1,\ldots,1, t^{(-j+1)/2},t^{(-j+3)/2},\ldots,t^{(j-1)/2}),\]
(which does not lie in our gauge group), we transform further to
$\alpha_{0\infty} = 0$, and
\begin{equation}\label{eqn:normalformextra}
\beta_{0\infty}=
 \vect{ P_0& q_0e_+\\
               \vect{r_0\\0}&-s +s_0e_+}
\end{equation}
\end{itemize}

These normal forms are unique up to the action of $Gl(k,\C)$, if in
addition one asks that the ``extra data'' vector $v$ be mapped to the
$(k+1)$-th basis vector in the normal form.
\end{lemma}

One now must create a Nahm complex from the data of the sheaves and the exact
sequences. On the interior of the first interval, we
use the matrix $B$ coming from
Diagram (\ref{resolution-diagram}) and set
\begin{equation}
V_{\infty0}= (\theta_0, \theta_\infty)\times H^0(\PP^1,Q_{\infty0}), \quad
\alpha_{\infty0} = 0,\quad \beta_{\infty0} = -B.
\end{equation}
In the same vein, we set

\begin{equation}\label{complex-bigger side}\begin{gathered}
V_{0\infty} = (\theta_\infty,2\pi+ \theta_0)\times H^0(\PP^1,Q_{0\infty}),  \quad
\alpha_{0\infty} = 0,\\
\beta_{0\infty} =
  \vect{-B&C_1e_+\\
                 \vect{-B'\\0}&-s+ C'_1e_+}=-M.
\end{gathered}\end{equation}
Note that the endomorphism is already in the normal form
(\ref{eqn:normalformextra}) given by Lemma \ref{lemma:normalform}.
There remains the gluing on the ends of the interval, which
involves introducing some form of monodromy, as we are on a circle.
The gluing will be mediated by the sheaves $P_0,P_\infty$.

The basic trick is that, for $n\leq m$, there are inclusions
\begin{align*}
f^0_{m,n}\colon H^0(\PP^1, P_0(n))&\rightarrow H^0(\PP^1, P_0(m)),\\
f^\infty_{m,n}\colon H^0(\PP^1, P_\infty(n))&\rightarrow
H^0(\PP^1, P_\infty(m)),
\end{align*}
as sections vanishing to an appropriate order at infinity. In
addition, there are maps arising from the exact sequences
(\ref{sheaf-sequences2})
\begin{align*}
n^0_\ell\colon H^0(\PP^1, P_0(\ell))&\rightarrow H^0(\PP^1,Q_{\infty0}),\\
n^\infty_\ell\colon H^0(\PP^1,P_\infty(\ell))&\rightarrow H^0(\PP^1,Q_{\infty0}),\\
m^0_\ell\colon H^0(\PP^1, P_0(\ell))&\rightarrow H^0(\PP^1,Q_{0\infty}),\\
m^\infty_\ell\colon H^0(\PP^1, P_\infty(\ell))&\rightarrow H^0(\PP^1,Q_{0\infty}).
\end{align*}
The maps $m^0_{-1}, n^0_{-j-1}, m^\infty_{j-1}, n^\infty_{-1}$ are
isomorphisms. Schematically,
\begin{diagram}
&&&H^0(Q_{\infty0})=\C^k\\
&\ruTo^{n^0_{-j-1}}_\cong&&&&\luTo>{n_{-1}^\infty}<\cong\\
H^0(P_0(-j-1))&&&&&&H^0(P_\infty(-1))\\
\dInto<{f^0_{-1,-j-1}}&&&&&&\dInto>{f^\infty_{j-1,-1}}\\
H^0(P_0(-1))        &&&&&&H^0(P_\infty(j-1))\\
&\rdTo<{m^0_{-1}}>\cong&&&&\ldTo>{m_{j-1}^\infty}<\cong\\
&&&H^0(Q_{0\infty})=\C^{k+j}.
\end{diagram}

At $\theta_0$, we define maps
\begin{align*}
\pi_0\colon V_{0\infty}(\theta_0)=H^0(\PP^1,Q_{0\infty})&\rightarrow V_{\infty0}(\theta_0)= H^0(\PP^1,Q_{\infty0})\\
i_0\colon V_{\infty0}(\theta_0)= H^0(\PP^1,Q_{\infty0})&\rightarrow
V_{0\infty}(\theta_0)=H^0(\PP^1,Q_{0\infty})
\end{align*}
by $i_0 = m^0_{-1}\circ f^0_{-1, -j-1}\circ (n^0_{-j-1})^{-1}$,
$\pi_0 = n^0_{-1}\circ (m^0_{-1})^{-1}$. The composition
$\pi_0\circ i_0$ is
the identity. Similarly, at $\theta_\infty$, we define maps
\begin{align*}
\pi_\infty\colon V_{0\infty}(\theta_\infty)=H^0(\PP^1,Q_{0\infty})&\rightarrow
V_{\infty0}(\theta_\infty)= H^0(\PP^1,Q_{\infty0})\\
i_\infty\colon V_{\infty0}(\theta_\infty)= H^0(\PP^1,Q_{\infty0})&\rightarrow
V_{0\infty}(\theta_\infty)=H^0(\PP^1,Q_{0\infty})
\end{align*}
by
$i_\infty = m^\infty_{j-1}\circ f^\infty_{j-1, -1}\circ (n^\infty_{-1})^{-1}$,
$\pi_\infty = n^\infty_{j-1}\circ (m^\infty_{j-1})^{-1}$. Again,
$\pi_\infty \circ i_\infty$ is the identity.

In the bases used in  Diagram (\ref{resolution-diagram}),
one has the block decompositions
\begin{align*}
i_0&=\vect{1\\0},\                           &\pi_0&=\vect{1&0},\\
i_\infty&=  \vect{A\\A'}=N\cdot\vect{1\\0}, &\pi_\infty&=\vect{1&0}\cdot N^{-1}.
\end{align*}

We set
\begin{equation*}
\vect{\tilde C_1\\{\tilde C}'_1}:= N^{-1}M^j\vect{C_2\\C'_2}\text{ and }
\tilde \beta_{0\infty} :=\vect{-B&\tilde C_1e_+\\
\vect{-D_2\\0}&-s+{\tilde C}'_1e_+},\end{equation*}
and then the monad equations imply
\begin{equation*}
\beta_{0\infty}N = N\tilde \beta_{0\infty}.\end{equation*}

We can then introduce $N$ as the parallel transport from $\theta_\infty$ to
$\theta_0+2\pi$ over the big side, as well
as introducing
the necessary poles. Indeed over the interval
$(\theta_\infty, 2\pi + \theta_0)$, we begin with the complex
$\alpha_{0\infty}, \beta_{0\infty}$ of Equation (\ref{complex-bigger side}),
and then
gauge it with a transformation $g$, given by choosing a smooth path $h(\theta)$
in $Gl(k+j,\C)$ equal to $N$ at $\theta_\infty+ 2\epsilon$ and
the identity at $2\pi + \theta_0-2\epsilon$, and setting
\[g(\theta)=\begin{cases} G(\theta-\theta_\infty)\circ N^{-1},
            & t\in (\theta_\infty, \theta_\infty+\epsilon),\\
     h(\theta),
    &  t\in (\theta_\infty+ 2\epsilon, 2\pi + \theta_0-2\epsilon),\\
    G(2\pi+\theta_0-\theta)^{-1},
        &  t\in (2\pi + \theta_0-\epsilon, 2\pi + \theta_0).
\end{cases}\]
We then smooth $g$ over the remaining small intervals, so that the result
is $C^\infty$. Applying $g$ to our Nahm complex over the interval,
we obtain
\begin{equation*}
(\gamma_{0\infty}, \delta_{0\infty}) = \Biggl( -\frac12\dot g g^{-1}, g \vect{-B&C_1e_+\\
                 -\vect{B'\\0}&-s+ C'_1e_+
   } g^{-1} \Biggr).
\end{equation*}
Under the gauge transformation $g$, the gluing maps become
\begin{equation*}
i_0= i_\infty=\vect{1\\0},\text{ and } \pi_0= \pi_\infty=\vect{1&0}.
\end{equation*}

The Nahm complex over the circle associated to the complex of sheaves is then
given by
\begin{gather*}\begin{aligned}
(\alpha_{\infty0}, \beta_{\infty0}) &= (0,B),\\
(\alpha_{0\infty},\beta_{0\infty}) &= (\gamma_{0\infty}, \delta_{0\infty}),
\end{aligned}\\
\begin{aligned}
i_0&= \vect{1\\0},
&i_\infty&=\vect{1\\0},\\
 \pi_0&=\vect{1&0},
&\pi_\infty&=\vect{1&0}.
\end{aligned}\end{gather*}

The ``extra data'' vectors are obtained from the
trivializations of $P_0$, $P_\infty$ at infinity.

Conversely, given a Nahm complex, we can recover the matrix data, and
hence the sheaves, by first gauging $\alpha_{\infty0}$ to zero, and
setting $B = \beta_{\infty0}$. Next gauging
$(\alpha_{0\infty},\beta_{0\infty})$ to their normal form near the poles,
as in the lemma above, so that $\alpha_{0\infty}= 0$ near the poles, gives
the matrix data $B', C_1, C'_1$ from the normal form near $\theta_0$, and
$D_2$ from the normal form near $\theta_\infty$. The gauge transformation
relating the two normal forms (that is the integration of the connection
$\alpha_{0\infty}$ between $\theta_0$  and  $\theta_\infty$) is the
matrix $N$ defined above; from it, one can recuperate $A,A', C_2,C'_2$.
Finally, setting $D_1 = e_+A'$ gives us the remaining data. The fact that
the matrix $N$ conjugates one normal form to the other then yields back the
monad equations.

For $j= 0$, the correspondence is much simpler. We build our bundles, and the
maps $i_0, \pi_0, i_\infty,\pi_\infty$ in the same way. The resolutions for
$Q_{\infty0}, Q_{0\infty}$ then give us matrices $B$, $\tilde B$, and it is
straight forward to see that $B - i_0\tilde B\pi_0$,
$B - i_\infty \tilde B \pi_0$ are of rank one. In the trivialization given
in the previous sections,
\begin{equation*}
          i_0 = 1, i_\infty = A.\end{equation*}
We then have, as above, the rank one jumps
\begin{equation*}
          B-\tilde B = C_1D_1 A^{-1}, \text{ and }B-A^{-1}\tilde B A = -
          A^{-1}C_2D_2.
\end{equation*}
One defines the Nahm complex by choosing a path $g(t)$ from the
identity to $A$, and setting
\begin{equation}
   \alpha_{\infty0} = 0, \beta_{\infty0} = B,\
   \alpha_{0\infty} = -\frac12\dot g g^{-1},\
   \beta_{\infty0} = gBg^{-1}.
\end{equation}

Again, the trivializations of $P_0, P_\infty$ at infinity give us the normalizations
of the decompositions of the jumps as a product of a column and a row.

Conversely,
from the Nahm complex, it is straightforward to extract the matrix
information, and so the sheaves.

The final correspondence which must be checked is the irreducibility
conditions. A Nahm complex over the circle is \emph{reducible} if there exists
a subbundle of each $V_*$, parallel for the $\alpha_*$ and invariant under the
$\beta_*$, mapping to each other by the gluing maps at the boundary point,
and proper on at least one interval.

Let us consider the three cases of reducibility for the complex of sheaves.
Because of the way the gluing maps, the connections $\alpha_{0\infty}$,
$\alpha_{\infty0}$ and the endomorphisms $\beta_{0\infty}$,
$\beta_{\infty0}$ are built for the cohomology of the $P_*$ and $Q_*$, we can
see that
\begin{itemize}
\item case 1 corresponds to the existence of a sub-line bundle of the $V_*$,
invariant and parallel,
\item case 2 corresponds to the existence of a co-rank 1 subbundle of the
$V_*$, invariant and parallel,
\end{itemize}

Summarizing:
\begin{theorem} \label{thm:SequencesComplexes}
Let $k\ge 1, j\ge 0$ be integers. There is an equivalence between

1) Exact sequences of sheaves
\begin{equation*}
\begin{matrix}
0&\rightarrow& \OO(j)&\rightarrow& P_0&\rightarrow& Q_{\infty0}&\rightarrow&0\cr
0&\rightarrow& \OO&\rightarrow& P_\infty&\rightarrow& Q_{\infty0}&\rightarrow&0\cr
0&\rightarrow& \OO&\rightarrow& P_0&\rightarrow& Q_{0\infty}&\rightarrow&0\cr
0&\rightarrow& \OO(-j)&\rightarrow& P_\infty&\rightarrow& Q_{0\infty}&\rightarrow&0
\end{matrix}
\end{equation*}
on $\PP^1$, with $Q_{0\infty}, Q_{\infty0}$ torsion sheaves of
length $k+j, k$ respectively, supported away from infinity, and with
$P_0, P_\infty$ trivialized over infinity on the line, and with
irreducible.

2) Irreducible Nahm complexes $\alpha, \beta$ on the circle, with rank $k$
over $(\theta_0,\theta_\infty)$, rank $k+j$ over
$(\theta_\infty, 2\pi + \theta_0)$, with the boundary conditions defined above,
modulo the action of the complex gauge group.
\end{theorem}

The last major step is to pass from Nahm complexes to solutions of Nahm's
equations. These equations are obtained by adding to the covariant
constancy condition
\begin{equation}\label{eqn:Nahmcomplex}
      \frac{d\beta}{dt} + [\alpha,\beta] = 0,
\end{equation}
the additional ``real'' equation
\begin{equation}\label{eqn:Nahmreal}
  \frac{d(\alpha+\alpha^*)}{dt} + [\alpha,\alpha^*] + [\beta, \beta^*] = 0.
\end{equation}
These equations are invariant under {\it unitary} gauge
transformations. One then has the theorem that orbits of irreducible
Nahm complexes under the action of the complex gauge group contain a
unique solution to Nahm's equations, up to  the action of the
unitary gauge group. The idea of the proof, due to Donaldson
\cite{DonaldsonSU2}, is to give a variational formulation to the
equations, and to show that each orbit contains a unique critical
point. The proof given in \cite[Sect. 2]{hurtubiseClassification} in
the context of $SU(N)$ monopoles extends verbatim to the case we
consider here, with one main difference, that of irreducibility.
In \cite[Sect. 2]{hurtubiseClassification}, the irreducibility is
automatic, because of the pole structure. Here, as we have seen
for the Nahm complexes, the irreducibility must be put in as a
supplementary condition. We note that, with the addition of a
metric structure, the notions of codimension one invariant
subbundle and dimension one invariant subbundle fuse as one
solves the variational problem, into a dimension one s
ubbundle, invariant under the $T_i$; in other words,
minimizing the energy takes one from a block upper triangular
form or block lower triangular form to a block diagonal form.

This last equivalence, combined with Theorems
\ref{thm:BundleMatrices}, \ref{thm:MatricesSequences}, and
\ref{thm:SequencesComplexes}, provides a proof of our main result,
Theorem \ref{MainTheorem}, which we now rewrite in the language we
absorbed throughout our journey.

\begin{theorem}
Let $k\ge 0, j\ge 0$ be integers. There is an equivalence between

1) Vector bundles $E$ of rank two on $\PP^1\times \PP^1$, with
   $c_1(E) = 0, c_2(E) =k$.
trivialized along $\PP^1\times\{\infty\}\cup\{\infty\}\times\PP^1$,
and with a based flag $\phi\colon\OO(-j)\hookrightarrow E$ of degree
$j$ along $\PP^1\times\{0\}$ (up to non-zero scalar multiple),
with the basing condition $\phi
(\infty) (\OO(-j)) = {\rm span}(0,1)$, and

2) Irreducible solutions to Nahm's equations $\alpha, \beta$ on the circle,
with rank $k$ over $(\theta_0,\theta_\infty)$, rank $k+j$ over
$(\theta_\infty, 2\pi + \theta_0)$, with the boundary conditions defined above,
modulo the action of the unitary gauge group.
\end{theorem}

The case $k=0$ is a version of a theorem of Donaldson \cite{DonaldsonSU2}, and
the other cases have been dealt with above.






\def\cprime{$'$}
\providecommand{\bysame}{\leavevmode\hbox to3em{\hrulefill}\thinspace}
\providecommand{\href}[2]{#2}

\where

\end{document}